\documentclass[12pt,english]{article}

\usepackage{hyperref}
\hypersetup{
    unicode      = false,
    pdftoolbar   = true,
    pdfmenubar   = true,
    pdffitwindow = true,
    pdfnewwindow = true,
    colorlinks   = true,
    linkcolor    = blue,
    citecolor    = blue,
    filecolor    = blue,
    urlcolor     = blue,
    breaklinks   = true
}

\usepackage[utf8]{inputenc}
\usepackage[T1]{fontenc}
\usepackage{amsthm,amssymb,amsmath}

\usepackage[english]{babel}
 \usepackage{lmodern}
\usepackage{graphicx}
\usepackage{fullpage}
\usepackage{caption}
\usepackage{version}
\usepackage{subcaption}
\usepackage{mathdots}

\usepackage{lipsum}
\usepackage{geometry}
\usepackage{systeme}
\usepackage{array,multirow,makecell}
\setcellgapes{1pt}
\makegapedcells
\usepackage[table]{xcolor}
\newcolumntype{R}[1]{>{\raggedleft\arraybackslash }b{#1}}
\newcolumntype{L}[1]{>{\raggedright\arraybackslash }b{#1}}
\newcolumntype{C}[1]{>{\centering\arraybackslash }b{#1}}
\usepackage{graphicx}
\usepackage{tikz}
\usetikzlibrary{positioning,decorations.pathreplacing,shapes}
\usepackage{pict2e}
\usepackage{amsfonts}
\usepackage{multirow}
\usepackage[normalem]{ulem}
\usepackage{epsfig}
\usepackage{soul}
\usepackage{colortbl}
\usepackage{enumitem}
\usepackage{optidef}
\usepackage{epic}
\usepackage{xspace}
\usepackage{pgf}
\usepackage{rotating}
\usepackage{array}
\usepackage{geometry}
\usepackage{times}
\usepackage{algorithm}
\usepackage{algpseudocode}
\usepackage[justification=centering]{caption}
\usepackage{nomencl}

\usetikzlibrary{patterns}
\usetikzlibrary{arrows}

\makeatletter
\tikzset{
        hatch distance/.store in=\hatchdistance,
        hatch distance=5pt,
        hatch thickness/.store in=\hatchthickness,
        hatch thickness=5pt
        }
\newcounter{subsubsubsection}[subsubsection]
\renewcommand\thesubsubsubsection{\@roman\c@subsubsubsection}
\newcommand\subsubsubsection{\@startsection{subsubsubsection}{4}{\z@}%
                                     {-3.25ex\@plus -1ex \@minus -.2ex}%
                                     {1.5ex \@plus .2ex}%
                                     {\normalfont\small\bfseries}}
\newcommand*\l@subsubsubsection{\@dottedtocline{3}{5.2em}{1em}}
\newcommand*{\subsubsubsectionmark}[1]{}

\makeatother

\def\m{\mathbf{m}}
\def\g{\mathbf{g}}

\def\x{\mathbf{x}}
\def\u{\mathbf{u}}
\def\y{\mathbf{y}}

\def\R{{\mathbb{R}}}
\def\N{{\mathbb{N}}}

\newcommand{\argmin}{\mathop{\mathrm{argmin}}}

\newtheorem{Def}{Definition}[section]
\newtheorem{Lem}[Def]{Lemma}
\newtheorem{Prop}[Def]{Proposition}
\newtheorem{Th}[Def]{Theorem}
\newtheorem{Cor}[Def]{Corollary}
\newtheorem{assumption}{Assumption}

\title{Sequential stochastic blackbox optimization with zeroth-order gradient estimators}

\author{
    \addtocounter{footnote}{1}
	\href{mailto:Charles.Audet@gerad.ca}{Charles Audet}
	\thanks{
		{GERAD}
		and D\'epartement de math\'ematiques et g\'enie industriel,
		\'Ecole Polytechnique de Montr\'eal,
		Montr\'eal, Qu\'ebec, Canada \newline
		Mail: charles.audet@gerad.ca
	}
	\and
	\href{mailto:jean.bigeon@ls2n.fr}{Jean Bigeon}
	\thanks{
	Nantes University, École Centrale Nantes, CNRS, LS2N, UMR 6004, F-44000 Nantes, France  \newline
	Mail: jean.bigeon@ls2n.fr
	}
	\and 
	\href{mailto:romain.couderc@grenoble-inp.fr}{Romain Couderc} \footnotemark[2] 
	\thanks{Univ. Grenoble Alpes, CNRS, Grenoble INP*, G-SCOP, 38000 Grenoble, France. \newline
	Mail: romain.couderc@grenoble-inp.fr \newline
    *Institute of Engineering Univ. Grenoble Alpes
    }
    \and
    \href{mailto:michael.kokkolaras@mcgill.ca}{Michael Kokkolaras} 
    \thanks{GERAD and Department of Mechanical Engineering, McGill University, Montreal, Canada \newline
    Mail: michael.kokkolaras@mcgill.ca}
}

\date{\today}

\begin{document}

\maketitle

\begin{abstract}
This work considers stochastic optimization problems in which the objective function values can only be computed by a blackbox corrupted by some random noise following an unknown distribution. 
The proposed method is based on sequential stochastic optimization (SSO): the original problem is decomposed into a sequence of subproblems. 
Each subproblem is solved using a zeroth-order version of a sign stochastic gradient descent with momentum algorithm (ZO-Signum) and with an increasingly fine precision. This decomposition allows a good exploration of the space while maintaining the efficiency of the algorithm once it gets close to the solution. 
Under Lipschitz continuity assumption on the blackbox, a convergence rate 
in expectation is derived for the ZO-Signum algorithm. 
Moreover, if the blackbox is smooth and convex or locally convex around its  minima, a convergence rate to an $\epsilon$-optimal point of the problem may be obtained for the SSO algorithm. 
Numerical experiments are conducted to compare the SSO algorithm with other state-of-the-art algorithms and to demonstrate its competitiveness.
\end{abstract}

\textbf{Mathematics Subject Classification:} 90C15, 90C56, 90C30, 90C90, 65K05

\section{Introduction}

The present work targets stochastic blackbox optimization problems of the form
\begin{equation}
    \min_{\x \in \R^n} \;  f(x) 
    \qquad \mbox{ where } \qquad
    f(\x) := \mathbb{E}_{\boldsymbol{\xi}} \left[ F(\x, \boldsymbol{\xi}) \right] ,
    \label{prob1}
\end{equation}
and in which 
    $F:\R^n \times \R^m \to \R$ is a blackbox \cite{AuHa2017}
    that takes two inputs:
    a vector of design variables $\x \in \R^n$ and
    a vector $\boldsymbol{\xi} \in \R^m$ that represents
    random uncertainties with an unknown distribution.
The function $F$ is called a {\em stochastic zeroth-order oracle} \cite{ghadimi2013stochastic}.
The objective function $f$ is obtained by taking the expectation of $F$
 over all possible values of the uncertainties $\boldsymbol{\xi}$.
The main applications belong to two different fields. 
The first is in a machine learning framework where the loss function’s gradient is unavailable or difficult to compute, for instance in optimizing neural network architecture \cite{real2017large}, design of adversarial attacks \cite{chen2019zo}, or game content generation \cite{volz2018evolving}. 
The second field is when the function $F$ is evaluated by means of a computational procedure~\cite{Kokko06}. In many cases, it depends on an uncertainty vector $\boldsymbol{\xi}$ due to environmental conditions, costs, or effects of repair actions that are unknown \cite{rockafellar2015risk}. Another source of uncertainty appears when the optimization is conducted at the early stages of the design process, where knowledge, information, and data is very limited.

\subsection{Related work}

Stochastic derivative-free optimization has been the subject of research for many years. Traditional derivative-free methods may be divided into two categories \cite{CoScVibook}: direct search and model-based methods. Algorithms corresponding to both  methods have been adapted to stochastic zeroth-order oracle. Examples include the stochastic Nelder-Mead algorithm \cite{chang2012stochastic} and the stochastic versions of the MADS algorithm \cite{audet2021stochastic, AudIhaLedTrib2016} for the direct search methods. For model-based methods, most work consider extensions of the trust region method \cite{chen2018stochastic, curtis2019stochastic, maggiar2018derivative}. A major shortcoming of these methods is their difficulty to  scale to large problems.

Recently, another class of methods, named zeroth-order (ZO) methods, has been attracting increasing attention. 
These methods use stochastic gradient estimators, which are based on the seminal work in \cite{kiefer1952stochastic, robbins1951stochastic} and have been extended in \cite{ghadimi2013stochastic, nesterov2017random, rubinstein_simulation_1981, spall_multivariate_1992}. 
These estimators have the appealing property of being able to estimate the gradient with only one or two function evaluations, regardless of the problem size. Zeroth-order methods take advantage of this property to extend first-order methods. For instance, the  well known first-order methods Conditional Gradient (CG), sign Stochastic Gradient Descent (signSGD) \cite{bernstein2018signsgd} and ADAptive Momentum (ADAM) \cite{kingma2014adam} have been extended to ZSCG \cite{balasubramanian2022zeroth}, ZO-signSGD \cite{liu2018signsgd} and ZO-adaMM \cite{chen2019zo}, respectively. More methods, not only based on first-order algorithms, have also emerged to solve regularized optimization problem \cite{cai2022zeroth}, for very high dimensional blackbox optimization problem \cite{cai2021zeroth} and for stochastic composition optimization problem \cite{ghadimi2020single}. 
Methods using second-order information based limited function queries have been developed~\cite{kim2021curvature}. Some methods handle situations where the optimizer has only access to a comparison oracle which indicates which of two points has the highest value \cite{cai2022one}. For an overview on ZO methods, readers may consult \cite{liu2020primer}.

\subsection{Motivation}

\begin{figure}[htb!]
     \centering
     \includegraphics[width = 0.5 \linewidth]{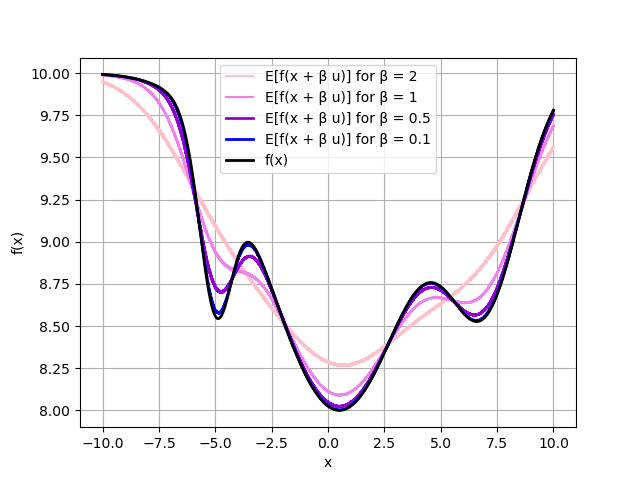}
     \caption{Curves of $f^\beta$ for $u \sim \mathcal{N}(0, 1)$ and different values of $\beta$.}
     \label{fig_expect}
 \end{figure}

 Formally, stochastic gradient estimators involve a smoothed functional $f^\beta$ (see Chapter 7.6 in \cite{rubinstein_simulation_1981}) which is a convolution product between $f$ and a kernel $h^\beta(\u)$
 \begin{equation}
     f^\beta(\x) := \int_{-\infty}^\infty h^\beta(\u)f(\x - \u) d \u = \int_{-\infty}^\infty h^\beta(\x - \u)f(\u) d \u.
 \end{equation}
 The kernel must fulfill a set of conditions [pp. 263, \cite{rubinstein_simulation_1981}]:
 \begin{enumerate}
 \label{condition}
     \item $h^\beta(\u) = \frac{1}{\beta^n} h(\frac{\u}{\beta})$ is a piecewise differentiable function;
     \item $\lim_{\beta \to 0} h^\beta(\u) = \delta(\u)$, where $\delta(v)$ is Dirac's delta function; 
     \item $\lim_{\beta \to 0} f^\beta(\x) = f(\x) $, if $\x$ is a point of continuity of $f$; \label{cond3}
     \item The kernel $h^\beta(\u)$ is probability density function (p.d.f.), that is $f^\beta(\x) = \mathbb{E}_{\mathbf{U} \sim h^\beta(\u)}[f(\x - \mathbf{U} )] = \mathbb{E}_{\mathbf{U} \sim h(\u)}[f(\x - \beta U)]$.
 \end{enumerate}
 Frequently used kernels include the Gaussian distribution and the uniform distribution on a unit ball. 
 Three properties concerning smoothed functional are worth noting. 
 First, the smoothed functional may be interpreted as a local weighted average of the function values in the neighbourhood of $\x$. 
 Condition \ref{cond3} implies that it is possible to obtain a solution arbitrarily close to a local minimum $f^*$. 
 Second, the smoothed functional is infinitely differentiable as a consequence of the convolution product, regardless of the degree of smoothness of $f$. 
 Moreover, according to the chosen kernel, stochastic gradient estimators may be calculated. 
 These estimators are unbiased estimators of $\nabla f^\beta$ and may be constructed on the basis of observations of $F(\x, \boldsymbol{\xi})$ alone. 
 Finally, the smoothed functional allows  convexification of the original function $f$. Previous studies \cite{rubinstein81, styblinski1990experiments} show that greater values of $\beta$ result in better convexification, as illustrated in Figure \ref{fig_expect}. Additionally, a larger $\beta$ leads to greater exploration of the space during the calculation of the gradient estimator. It has also been demonstrated in \cite{liu2018zeroth} that if the smoothing parameter is too small, the difference in function values cannot be used to accurately represent the function differential, particularly when the noise level is significant.

 Although the two first properties of the smoothed functional are exploited by ZO methods, the last property has not been utilized since the work in \cite{styblinski1990experiments}. This may be because the convexification phenomenon becomes insignificant when dealing with high-dimensional problems \footnote{Note that a blackbox optimization problem with dimensions ranging from $100$ to $1000$ may be considered large, while problems with  $n \geq 10000$ may be considered very large.}. However, for problems of relatively small size ($n \simeq 10$), this property can be useful.  The authors of ~\cite{styblinski1990experiments} use an iterative algorithm to minimize the sequence of subproblems
 \begin{equation}
     \min_{\x \in \R^n} f^{\beta^i}(\x),
     \label{subproblems}
 \end{equation}
where $\beta^i$ belongs to a finite prescaled sequence of scalars. This approach is limited because the sequence $\beta^i$ does not necessarily converge to $0$ and the number of iterations to go from subproblem $i$ to $i+1$ is arbitrarily fixed a priori. Furthermore, neither a convergence proof nor a convergence rate are provided for the algorithm. Finally, although promising, numerical results are only presented for analytical test problems. These shortcomings motivate the research presented here. 

\subsection{Contributions}
The main contributions of this paper can be summarized as follows. 
\begin{itemize}
    \item A sequential stochastic (SSO) optimization algorithm is developed to solve the sequence of subproblems in Equation~(\ref{subproblems}). In the inner loop,  a subproblem is solved according to the zeroth-order version of the Signum algorithm \cite{bernstein2018signsgd}. The stopping criterion is based on the norm of the momentum, which must be below a certain threshold. In the outer loop, the sequence of $\beta^i$ is proportional to the threshold needed to consider a subproblem solved and is driven to $0$. Therefore, the smaller the value of $\beta^i$ is (and thus better is the approximation given by $f^{\beta^i}$), the larger the computational budget granted to the resolution of the subproblem. 
    
    \item A theoretical analysis of this algorithm is conducted. First, the norm of the momentum is proved to converge to $0$ in expectation, with a convergence rate that depends on the step sizes. Then, the convergence rate in expectation of the ZO-Signum algorithm to a stationary point of $f^\beta$ is derived under Lipschitz continuity of the function $F$. Finally, if the function $F$ is smooth, and $f^{\beta}$ is convex or become convex around its local minima, a convergence rate to an $\epsilon$-optimal point is derived for the SSO algorithm.
    
    \item Numerical experiments are conducted to evaluate the performance of the proposed algorithm in two applications. First,  a comparison is made with traditional derivative-free algorithms on the optimization of the storage cost of a solar thermal power plant model, which is a low-dimensional problem. Second, a comparison is made with other ZO algorithms in order to generate blackbox adversarial attacks, which are large size problems.
\end{itemize}

The remainder of this paper is organized as follows. 
In Section~\ref{sec-GaussianGradientEstimator}, the main assumptions and the Gaussian gradient estimator are described. 
In Section~\ref{sec-SSO}, the sequential optimization algorithm is presented, and its convergence properties are studied in Section~\ref{sec-conv_analysis}. 
Section~\ref{sec-num_res} presents numerical results, and Section~\ref{sec-conclusion} draws conclusions and discusses future work. 

\section{Gaussian gradient estimator}
\label{sec-GaussianGradientEstimator}

The assumptions concerning the stochastic blackbox function $F$ are as follows.
\begin{assumption}
\label{assum1}
Let $(\Omega, \mathcal{F}, \mathbb{P})$ be a probability space. 
\begin{enumerate}[label=\alph*.]
    \item The function satisfies $F(\cdot, \boldsymbol{\xi}) \in L^1(\Omega, \mathcal{F}, \mathbb{P})$ and $f(\x) := \mathbb{E}_{\boldsymbol{\xi}}[F(\x, \boldsymbol{\xi})] $, for all $\x \in \R^n$.
    \item $F(\cdot, \boldsymbol{\xi})$ is Lipschitz continuous for any $\boldsymbol{\xi}$, with  constant $L_0(F)  > 0$.
\end{enumerate}
\end{assumption}
Assumption \ref{assum1}.a implies that the expectation of $F(\x, \boldsymbol{\xi})$ with respect to $\boldsymbol{\xi}$ is well-defined on  $\R^n$ and that the estimator $F(\x, \boldsymbol{\xi})$ is unbiased. 
Assumption \ref{assum1}.b is commonly used to ensure convergence and to bound the variance of the stochastic zeroth-order oracle. It is worth noticing that no assumption is made on the differentiability of the objective function $f$ or of its estimate $F$ with respect to $\x$, contrary to most work on zeroth-order methods.

Under Assumption \ref{assum1}, a smooth approximation of the function $f$ may be constructed by its convolution with a Gaussian random vector. Let $\u$ be an n-dimensional
standard Gaussian random vector and $\beta >0$ be the smoothing parameter. Then, a
smooth approximation of $f$ is defined as
\begin{equation}
    f^{\beta}(\x) := \frac{1}{(2  \pi)^{\frac{n}{2}}}\int f(\x +\beta \u) e^{-\frac{||u||^2}{2}} d\u = \mathbb{E}_\u[f(\x+\beta \u)].
\end{equation}

 This estimator has been studied in the literature (especially in \cite{nesterov2017random}) and benefits of several appealing properties. The properties used in this work are summarized in the following Lemma.

 \begin{Lem}
     Under Assumption \ref{assum1}, the following statements hold for any integrable function $f:\R^n \to \R$ and its approximation $f^{\beta}$ parameterized by $\beta > 0$.
     \begin{enumerate}
         \item $f^\beta$ is infinitely differentiable: $f^\beta \in \mathcal{C}^\infty$.
         \item A one-sided unbiased estimator of $\nabla f^\beta$ is 
    \begin{equation}
        \Tilde{\nabla} f^\beta(\x) := \frac{\u (f(\x + \beta \u) -f(\x))}{\beta}.
        \label{estim1}
    \end{equation}
    \item Let $\beta^2 \geq \beta^1 \geq 0$, then $\forall \x \in \mathbb{R}^n$
    \begin{equation*}
        ||\nabla f^{\beta^1}(\x) - \nabla f^{\beta^2}(\x)|| \leq L_1(f^{\beta^1})  (\beta^2-\beta^1)  (n+3)^{\frac{3}{2}}.
    \end{equation*}
    Moreover, for $\beta >0$, then $f^\beta$ is $L_1(f^\beta)$-smooth, i.e, $f^\beta \in \mathcal{C}^{1+}$ with $L_1(f^\beta) = \frac{2 \sqrt{n}}{\beta}L_0(F)$.
    \item If $f$ is convex, then $f^\beta$ is also convex.
     \end{enumerate}
     \label{lem_trunc_gaus}
     
 \end{Lem}
 \begin{proof}
     1. It is a consequence of the convolution product between an integrable function and an infinitely differentiable kernel.\\

     2.  See \cite[Equation (22)]{nesterov2017random}. \\

     3. If $\u \sim \mathcal{N}(\mathbf{0}, \mathbf{I})$, let define for all $\x \in \R^n$
\begin{align*}
    g(\x) = f^{\beta^1}(\x) = \mathbb{E}_\u[f(\x + \beta^1 \u)].
\end{align*}
Let $\mu = \beta^2 - \beta^1 \geq 0$, it follows that for all $\x \in \R^n$
\begin{equation*}
    g^\mu(\x) = \mathbb{E}_\u[g(\x + \mu \u)] = \mathbb{E}_\u[ f^{\beta^1}(\x + \mu \u)] = \mathbb{E}_\u[ f(\x + \mu \u + \beta^1 \u)] = \mathbb{E}_\u[ f(\x + \beta^2 \u)] = f^{\beta^2}(\x).
\end{equation*}
Then, since by \cite[Lemma 2]{nesterov2017random} under Assumption \ref{assum1}, $f^{\beta^1}$ is Lipschitz continuously differentiable, the \cite[Lemma 3]{nesterov2017random} may be applied to the function $g$ and it follows that
\begin{align*}
    ||\nabla f^{\beta^1}(\x) - \nabla f^{\beta^2}(\x)|| =  || \nabla g(\x) - \nabla g^{\mu}(\x) || 
    \leq  L_1(f^{\beta^1}) \mu (n+3)^{\frac{3}{2}} =  L_1(f^{\beta^1}) (\beta^2 - \beta^1) (n+3)^{\frac{3}{2}}.\\
\end{align*}

    4. See \cite[page 5]{nesterov2017random}.
 \end{proof}
The estimator obtained in Equation (\ref{estim1}) may be adapted to the stochastic zeroth-order oracle $F$. For instance, a one-sided (mini-batch) estimator of the noised function $F$ is
\begin{equation}
        \Tilde{\nabla} f^\beta(\x, \boldsymbol{\xi}) = \frac{1}{q} \sum_{j = 1}^q \frac{\u^j (F(\x + \beta \u^j, \boldsymbol{\xi}^j) - F(\x, \boldsymbol{\xi}^j))}{\beta},
        \label{estim4}
    \end{equation}
where $(u^j)_{j=1}^q$ and $(\boldsymbol{\xi}^j)_{j=1}^q$ are  $q$ Gaussian random direction vectors and their associated $q$ estimates values of the function $F$. This is still an unbiased estimator of $\nabla f^\beta$ because 
\begin{equation}
    \mathbb{E}_{\u, \boldsymbol{\xi}}[\Tilde{\nabla} f^\beta(\x, \boldsymbol{\xi})] =  \mathbb{E}_{\u}[ \mathbb{E}_{\boldsymbol{\xi}}[\Tilde{\nabla} f^\beta(\x, \boldsymbol{\xi}) | \u] ] = \nabla f^\beta(\x).
    \label{unbiased}
\end{equation}

The result of Lemma \ref{lem_trunc_gaus}.3 is essential to understand why solving a sequence of optimization problems defined in Equation (\ref{subproblems}) may be efficient, while it might seem 
counterproductive at first sight. Below are examples of the advantages of treating the problem with sequential smoothed function optimization.
\begin{itemize}
    \item The subproblems are approximations of the original problem and it is not necessary to solve them exactly. Thus, an appropriate procedure for solving these problems with increasingly fine precision can be used. Moreover, as seen in Lemma \ref{lem_trunc_gaus}.3, the norm of the gradient obtained in a subproblem is close to the one of the following subproblem. The computational effort to find a solution to the second subproblem from the solution of the first should therefore not be important.  
    \item  The information collected during the optimization process of a subproblem may be reused in the subsequent subproblems since they are similar. 
    \item A specific interest in the case of smoothed functional is the ability of using a larger value of $\beta$ during the solving of the first subproblems. It allows for a better exploration of the space and convexification phenomenon of the function (see Figure \ref{fig_expect}). Moreover, the new step size may be used for each subproblem, it allows increasing the step size momentarily, in the hope of having a greater chance of escaping a local minimum.
\end{itemize}

\section{A Sequential  Stochastic Optimization (SSO) algorithm }
\label{sec-SSO}

Section 3.1 presents a zeroth-order version of the Signum algorithm \cite{bernstein2018signsgd} to solve Subproblem (\ref{subproblems}) for a given $\beta^i$ and Section 3.2 presents the complete algorithm used to solve the sequential optimization problem.

\subsection{The Zeroth-Order Signum algorithm}

\begin{algorithm}[htb!]
	\caption{ZO-Signum (ZOS) algorithm  to solve subproblem $i \in \mathbb{N} $} 
	\begin{algorithmic}[1]
        \State{\textbf{Input:} $\x^{i, 0}, \m^{i,0}, \beta^i, s_1^{i,0}, s_2^{i,0}$, $L$, $q$, $M$}
        \State{Set $k = 0$ }
        \State{Define stepsize sequences $s_1^{i,k} = \frac{s_1^{i,0}}{(k+1)^{\alpha_1}}$ and $s_2^{i,k} =  \frac{s_2^{i,0}}{(k+1)^{\alpha_2}}$}
        \While{ $||\m^{i, k}|| > \frac{L \beta^i }{4 \beta^0}$ or $k \leq M$}
        \State{Draw $q$ samples $\u^k$ from the Gaussian distribution $\mathcal{N}(\mathbf{0}, I)$}
        \State{Calculate the average of the $q$ Gaussian estimate $\Tilde{\nabla} f^{\beta^i}(\x^{i,k}, \boldsymbol{\xi}^{i,k})$ from Equation (\ref{estim4})}
        \State{Update:}
        \begin{align}
            \label{update_grad}
            &\m^{i, k+1} = s_2^{i,k} \Tilde{\nabla} f^{\beta^i}(\x^{i,k}, \boldsymbol{\xi}^{k}) + (1-s_2^{i,k})\m^{i,k}\\
            \label{update_x}
            &x_j^{i, k+1} = x_j^{i,k} - s_1^k \text{sign}(m_j^{i,k+1})  \; \forall j \in [1, n] 
        \end{align}
 
        \State{$k \leftarrow k +1$}
        \EndWhile
    \State{Return $\m^{i, k}$ and $\x^{i,k}$}
	\end{algorithmic} 
	\label{algo1}
\end{algorithm}

A zeroth-order version of the Signum algorithm (Algorithm 2 of \cite{bernstein2018signsgd}) is used to solve the subproblems. 
The Signum algorithm is a momentum version of the sign-SGD algorithm. In  \cite{liu2018signsgd}, the authors extended the original sign-SGD algorithm to a zeroth-order version of this algorithm. However, a zeroth-order version of Signum is not studied in the work of \cite{liu2018signsgd}. As the Signum algorithm has been shown to be competitive with the ADAM algorithm \cite{bernstein2018signsgd}, a zeroth-order version of this algorithm seems interesting to consider. For completeness, the versions of the signSGD and the Signum as they originally appeared in \cite{bernstein2018signsgd} are given in Appendix \ref{appendix_C}. There is an important difference between the original Signum algorithm and its zeroth-order version presented in Algorithm \ref{algo1}. 
Indeed, while the step size of the momentum $1-s_2$ is kept constant in the work of \cite{bernstein2018signsgd}, it is driven to $0$ in our work. 
\begin{algorithm}[htb!]
	\caption{Sequential Stochastic Optimization (SSO) algorithm} 
	\begin{algorithmic}[1]
        \State{\textbf{Initialization:}}
		\State{Set $\x^{0,0} \in \R^{n}$, $\beta^0 > 0$ and $N$ the maximum number of function calls for the search step}
        \State{Set $q$ the number of gradient estimates at each iteration of ZO-Signum algorithm}
        \State{Set $M$ the minimum number of iterations made by the ZO-Signum algorithm on a subproblem}
        \State{$\mathcal{C}$ the cache containing all the evaluated points} 
        \State{Set $\m^{0, 0}= \Tilde{\nabla} f^{\beta^0}(\x^{0, 0}, \boldsymbol{\xi}^0)$} and $L = + \infty$
        \State{Set $s_1^{0, 0} > 0$ and  $s_2^{0, 0} > 0$  }
        \State{Set $i = 0$}
        \State{\textbf{Search step (optional):}}
        \While{$M(i+1)q \leq  N$:}
            \State{Solve subproblem $i$ with Algorithm \ref{algo1}:}
            \begin{align*}
                &\m^{i+1, 0} = ZOS(\x^{i, 0}, \m^{i, 0},  \beta^i, s_1^{i, 0}, s_2^{i,  0}, L, q, M) \\
                &\x^{i+1, 0} \in \argmin_{\x \in \mathcal{C}} F(\x, \boldsymbol{\xi})
            \end{align*}
            \State{Update $\beta^i$, $s_1^{i,0}$ and $s_2^{i,0}$ as in step 18}
        \EndWhile
        \State{$ L = ||\m^{0,0}||$}
        \State{\textbf{Local step:}}
        \While{ $\beta^i > \epsilon$}
            \State{Solve subproblem $i$ with Algorithm \ref{algo1}:}
            \begin{equation*}
                \m^{i+1, 0}, \x^{i+1, 0} = ZOS(\x^{i, 0}, \m^{i, 0},  \beta^i, s_1^{i, 0}, s_2^{i,  0}, L, q, M)
            \end{equation*}

        \State{ Update:}
        \begin{align*}
            &\beta^i = \frac{\beta^0}{(i+1)^2}, 
            s_1^{i, 0} = \frac{s_1^{0, 0}}{(i  +1)^{\frac{3}{2}}}, \quad s_2^{i, 0} = \frac{s_2^{0,0}}{i+1}\\
            &i \leftarrow i + 1
        \end{align*}
            \EndWhile   
    \State{Return $\x^i$}
	\end{algorithmic} 
	\label{algo2}
\end{algorithm} 
This leads to two consequences. First, the variance is reduced since the gradient is averaged on a longer time horizon, without using mini-batch sampling. Second, as it has been demonstrated in other stochastic approximation works \cite[section 3.3]{bhatnagar_stochastic_2013}, \cite{ruszczynski1983stochastic}, with carefully chosen step sizes the norm of the momentum goes to $0$ with probability one. In the ZO-Signum algorithm, the norm of the momentum is thus used as a stopping criterion.

\subsection{The SSO algorithm}

The optimization of the subproblems sequence described in Equation (\ref{subproblems}) is driven by the sequential stochastic (SSO) algorithm presented in Algorithm \ref{algo2}. The value of $\beta$ plays a critical role, as it serves as both the smoothing parameter and the stopping criterion for Algorithms \ref{algo1} and \ref{algo2}. 
Algorithm \ref{algo2} is inspired by the MADS algorithm \cite{AuDe2006} as it is based on two steps: a search step and a local step. The search step is optional and may consist in any heuristics and is required only on problem with relatively small dimensions. In Algorithm \ref{algo2}, an example of a search is given which consists of updating $\x$ after $M$ iterations of the ZO-Signum algorithm with the best known $\x$ found so far. The local step is then used: Algorithm \ref{algo1} is launched on each subproblem $i$ with specific values of $\beta^i$ and step-size sequences. Once Algorithm \ref{algo1} meets the stopping criterion (which depends on the value of $\beta^i$), the value of $\beta^i$ and the initial step-sizes $s_1^{i,0}$ and $s_2^{i,0}$ are reduced, and the algorithm proceeds to the next subproblem. The convergence is guaranteed by the local step, since the search step is run only a finite number of times.

It is worth noting that the decrease rate of $\beta^i$ is chosen so that the difference between subproblems $i$ and $i+1$ is not significant. Therefore, the information collected in subproblem $i$, through the momentum vector $\m$, can be used in subproblem $i+1$. Furthermore, the initial step-sizes $s_1^{i,0}$ and $s_2^{i,0}$ decrease with each iteration, allowing us to focus our efforts quickly towards a local optimum when $s_1^{0,0}$ and $\beta^0$ are chosen to be relatively large. 

\section{Convergence analysis}
\label{sec-conv_analysis}

The convergence analysis is conducted in two steps : first the convergence in expectation is derived for Algorithm \ref{algo1}  and then the convergence for Algorithm \ref{algo2} is derived. 

\begin{table}[htb!]
    \centering
    \small
        \caption{Workflow of Lemmas/Propositions/Theorems for the ZO-Signum convergence analysis}
        \label{tab_conv_signum}
    \begin{tabular}{|c|c|c|c|c|}
        \hline
          Assumptions on $F$ & Preliminary results  & Intermediate results & Main results & When  $f^\beta$ is convex \\
         \hline
         & & \multirow{1}{*}{Proposition \ref{prop1}}& &
         \\ \cline{2-3}
          \multirow{2}{*}{Assumption \ref{assum1}} & \multirow{1}{*}{Lemma \ref{Lem_decomp}}  & \multirow{3}{*}{Proposition \ref{prop2}} & &\\
           \cline{2-2} 
          \multirow{2}{*}{which implies}& \multirow{1}{*}{Lemma \ref{lem_var}}& & \multirow{2}{*}{Theorem \ref{conv_expect}} & \multirow{2}{*}{Theorem \ref{th_convex}}\\
          \cline{2-2} 
         \multirow{2}{*}{$L_1(f^{\beta^i}) = \frac{2 \sqrt{n} L_0(F)}{\beta^i}$} &  \multirow{1}{*}{Lemma \ref{Lem_utile}} && \multirow{2}{*}{Corollary\ref{cor1}} & \\
         \cline{2-3} 
         & \multirow{1}{*}{Lemma \ref{Lem_reformulate}}  & \multirow{2}{*}{Proposition \ref{prop3}} & & \\
         \cline{2-2} 
         & \multirow{1}{*}{Lemma \ref{Lem_utile_bis}}  & & &\\
         \hline
    \end{tabular}
\end{table}

\subsection{Convergence of the ZOS algorithm}
The analysis of Algorithm \ref{algo1} follows the general methodology given in \cite[Appendix E]{bernstein2018signsgd}. In the following subsection, the main result  in \cite{bernstein2018signsgd} is recalled for completeness. The next subsections are devoted to bound  the variance and bias terms when  $\lim_{k \to \infty} s_2^{i,k} = 0$ . Finally, these results are used to obtain the convergence rate in expectation of Algorithm \ref{algo1} in the non-convex and convex case. The last subsection is devoted to a theoretical comparison with other ZO methods of the literature. The subproblem index $i$ is kept constant throughout this section. In order to improve the readability of the convergence analysis of the ZO-Signum algorithm, a hierarchical workflow of the different theoretical results is presented in Table \ref{tab_conv_signum}. The main results are Theorem~\ref{conv_expect} and its corollary for the nonconvex case, and~Theorem \ref{th_convex} for the convex case.

\subsubsection{Preliminary result \cite{bernstein2018signsgd}}
The following proposition uses the Lipschitz continuity of the function $f^{\beta^i}$ (proved in Lemma \ref{lem_trunc_gaus}) to bound the gradient at the $k$th iteration.
\begin{Prop}[\cite{bernstein2018signsgd}]
\label{prop1}
For the subproblem $i \in \N$, under Assumption \ref{assum1} and in the setting of Algorithm \ref{algo1}, we have
\begin{align}
\begin{split}
        s_1^{i,k} \mathbb{E}[||\nabla f^{\beta^i} (\x^{i,k}) ||_1] \leq &\mathbb{E}[f^{\beta^i}(\x^{i,k}) - f^{\beta^i}(\x^{i,k+1})]  + \frac{n L_1(f^{\beta^i})}{2} (s_1^{i,k})^2 \\
    &+ 2 s_1^{i,k} \underbrace{\mathbb{E}[||\Bar{\m}^{i,k+1}  - \nabla f^{\beta^i}(\x^{i,k})||_1]}_{\text{bias}} + 2s_1^{i,k}  \sqrt{n} \sqrt{\underbrace{\mathbb{E}[||\m^{i,k+1} - \Bar{\m}^{i,k+1}||_2^2]}_{\text{variance}}} 
\end{split}
\label{ineq_first}
\end{align}
where $\Bar{m}_j^{i,k+1}$ is defined recursively as $\Bar{m}_j^{i,k+1} = s_2^{i,k} \nabla f^{\beta^i}(\x^{i,k}) + (1-s_2^{i,k})\Bar{m}_j^{i,k} $.
\end{Prop}
\begin{proof}
    See Appendix \ref{appendix_A}.
\end{proof}
Now, it remains to bound the three terms on the right side of Inequality~(\ref{ineq_first}).
\subsubsection{Bound on the variance term}

The three following lemmas are consecrated to bound the variance term. Unlike the work reported in~\cite{bernstein2018signsgd}, the variance reduction is conducted by driving the step size of the momentum to $0$. It avoids to sample an increasing number of stochastic gradients at each iteration, which may be problematic, as noted in \cite{liu2018signsgd}. To achieve this, the variance term is first decomposed in terms of expectation of the squared norm of the stochastic gradient estimators $\Tilde{g}$.

\begin{Lem}
    For the subproblem $i \in \N$, let $k \in  \N$ and $j \in[1, n]$, we have 
    \begin{align*}
        \mathbb{E}[||\m^{i,k+1} - \Bar{\m}^{i,k+1}||^2]  &\leq (s_2^{i,k})^2 \mathbb{E}[ ||\Tilde{\g}^{i,k}||^2 ] + \sum_{r = 0}^{k-1} (s_2^{i,r})^2  \prod_{t = r}^{k-1} (1-s_2^{i,t+1})^2 \mathbb{E}[||\Tilde{\g}^{i,r}||^2]\\
        &+ \prod_{t = 0}^{k} (1-s_2^{i,t})^2 \mathbb{E}[||\Tilde{\g}^{i,0}||^2],
    \end{align*}
    where $\Tilde{g}_j^{i,r} = \Tilde{\nabla} f^{\beta^i}(\x^{i,r}, \boldsymbol{\xi}^r), \forall r \in [0, k]$ is defined in Equation (\ref{estim4}) and the norm is $||\cdot||_2$.
    \label{Lem_decomp}
\end{Lem}
\begin{proof}
    Let $k \in \N$, by definition of $\m^{i,k}$ and $\Bar{\m}^{i,k}$, it follows that
    \begin{align*}
        ||\m^{i,k+1} - \Bar{\m}^{i,k+1}||^2&= (s_2^{i,k})^2||\Tilde{\g}^{i,k} - \nabla f^{\beta^i}(\x^{i,k})||^2 + (1-s_2^{i,k})^2||\m^{i,k} - \Bar{\m}^{i,k}||^2\\
        &+2s_2^{i,k}(1-s_2^{i,k})( \Tilde{\g}^{i,k} - \nabla f^{\beta^i}(\x^{i,k}))^T (\m^{i,k} - \Bar{\m}^{i,k}).      
    \end{align*}
    The expectation of this expression is
    \begin{align}
         \mathbb{E}[||\m^{i,k+1} - \Bar{\m}^{i,k+1}||^2]&= (s_2^{i,k})^2\mathbb{E}[||\Tilde{\g}^{i,k} - \nabla f^{\beta^i}(\x^{i,k})||^2] + (1-s_2^{i,k})^2\mathbb{E}[||\m^{i,k} - \Bar{\m}^{i,k}||^2]\\
        &+2s_2^{i,k}(1-s_2^{i,k})\mathbb{E}[(\Tilde{\g}^{i,k} - \nabla f^{\beta^i}(x^{i,k}))^T (\m^{i,k} - \Bar{\m}^{i,k})].
        \label{exp}
    \end{align}
    Now, introducing the associated sigma field of the process $\mathcal{F}^{i,k} = \sigma(\x^{j,t}, \m^{j,t}, \Bar{\m}^{j,t}; j \leq i, t \leq k) $ by the law of total expectation, it follows that
    \begin{align*}
        \mathbb{E}[(\Tilde{\g}^{i,k} - \nabla f^{\beta^i}(x^{i,k}))^T (\m^{i,k} - \Bar{\m}^{i,k})] &= \mathbb{E}[ \mathbb{E}[(\Tilde{\g}^{i,k} - \nabla f^{\beta^i}(\x^{i,k}))^T (\m^{i,k} - \Bar{\m}^{i,k})| \mathcal{F}^{i,k}]] \\
        &= \mathbb{E}[ (\mathbb{E}[\Tilde{\g}^{i,k}|\mathcal{F}^{i,k}] - \nabla f^{\beta^i}(x^{i,k}))^T (\m^{i,k} - \Bar{\m}^{i,k})]\\
        &= 0,
    \end{align*}
    where the second equality holds because $\m^{i,k}, \Bar{\m}^{i,k}$ and $\nabla f^{\beta^i}(\x^{i,k})$ are fixed conditioned on $\mathcal{F}^{i,k}$ and because $\mathbb{E}[\Tilde{\g}^{i,k}| \x^{i,k}] = \nabla f^{\beta^i}(\x^{i,k}) $ as $\Tilde{\g}^{i,k}$ is an unbiased estimator of the gradient by Equation (\ref{unbiased}). By substituting this result in (\ref{exp}), it follows that
    \begin{equation*}
        \mathbb{E}[||\m^{i,k+1} - \Bar{\m}^{i,k+1}||^2]= (s_2^{i,k})^2\mathbb{E}[||\Tilde{\g}^{i,k} - \nabla f^{\beta^i}(\x^{i,k})||^2] + (1-s_2^{i,k})^2\mathbb{E}[||\m^{i,k} - \Bar{\m}^{i,k}||^2].
    \end{equation*}
    By repeating this process iteratively, we obtain 
    \begin{align}
    \begin{split}
        \mathbb{E}[||\m^{i,k+1} - \Bar{\m}^{i,k+1}||^2]  &= (s_2^{i,k})^2 \mathbb{E}[||\Tilde{\g}^{i,k} - \nabla f^{\beta^i}(\x^{i,k})||^2] \\
        &+ \sum_{r = 0}^{k-1} (s_2^{i,r})^2  \prod_{t = r}^{k-1} (1-s_2^{i,t+1})^2\mathbb{E}[||\Tilde{\g}^{i,r} - \nabla f^{\beta^i}(\x^{i,r})||^2]\\
        &+ \prod_{t = 0}^{k} (1-s_2^{i,t})^2 \mathbb{E}[||\Tilde{\g}^{i,0} - \nabla f^{\beta^i}(x^{i,0})||^2].
        \end{split}
        \label{res_int}
    \end{align}
    Finally, by observing that $\forall r \in [0,k], \mathbb{E}[\Tilde{\g}^{i,r}| \x^{i,r}] = \nabla f^{\beta^i}(x^{i,r})$ and  by the law of total expectation, we obtain
    \begin{align*}
              \mathbb{E}[||\Tilde{\g}^{i,r} - \nabla f^{\beta^i}(x^{i,r})||^2]  &= \mathbb{E}[ ||\Tilde{\g}^{i,r} - \mathbb{E}[\Tilde{\g}^{i,r}|\x^{i,r}] ||^2 ] \\
              &=  \mathbb{E}[|| \Tilde{\g}^{i,r} ||^2] - \mathbb{E}[||\nabla f^{\beta^i}(x^{i,r}) ||^2]\\
              &\leq \mathbb{E}[|| \Tilde{\g}^{i,r} ||^2].
    \end{align*}
    Introducing this inequality in Equation (\ref{res_int}) completes the proof.
\end{proof}
Second, the expectation of the squared norm of the stochastic gradient estimators are bounded by a constant depending quadratically on the dimension. 
\begin{Lem}
    Let $i \in \N$, $r \in [0,k]$, $j \in [1,n]$, then under Assumption \ref{assum1}, we have
    \begin{equation*}
        \mathbb{E}[|| \Tilde{\g}^{i,r} ||^2] \leq L_0(F)^2  (n + 4)^2
    \end{equation*}
    where $L_0(F)$ is the Lipschitz constant of $F$.
    \label{lem_var}
\end{Lem}
\begin{proof}
    By Equation (\ref{estim4}) with $q = 1$, it follows that  
    \begin{align*}
         \mathbb{E}[|| \Tilde{\g}^{i,r} ||^2] &= \mathbb{E}\left[\frac{||u||^2}{(\beta^i)^2} \left( F(\x^{i,r} + \beta^i \u, \boldsymbol{\xi}) - F(\x^{i,r}, \boldsymbol{\xi}) \right)^2 \right] \\
         &\leq  L_0(F)^2 \mathbb{E}[||u||^4] \\
         &\leq L_0(F)^2 (n+4)^2
    \end{align*}
   where the first inequality follows from Assumption \ref{assum1}.b and the second by \cite[Lemma 1]{nesterov2017random}.
\end{proof}

Finally, a technical lemma bounds the second term of the decomposition of the Lemma \ref{Lem_decomp} by a decreasing sequence. It achieves the same rate of convergence as in \cite{bernstein2018signsgd} without sampling any stochastic gradient.
\begin{Lem}
    For the subproblem $i \in \N$, let $s_2^{i,k}$ defined such that $s_2^{i,k} = \frac{s_2^{i,0}}{(k+1)^{\alpha_2}}$ with  $\alpha_2 \in (0, 1)$ and $s_2^{i,0} \in (0,1)$, then for $k$ such that 
    \begin{equation}
        \frac{k}{(k+1)^{\alpha_2}} \geq \frac{\ln(s_2^{i,0}) + (1+\alpha_2) \ln(k)}{s_2^{i,0}}
        \label{cond1}
    \end{equation}
    the following inequality holds
    \begin{equation}
        \sum_{r = 0}^{k-1} (s_2^{i,r})^2  \prod_{t = r}^{k-1} (1-s_2^{i,t+1})^2 \leq  \frac{9 s_2^{i,0} }{k^{\alpha_2}}.
    \end{equation}
    \label{Lem_utile}
\end{Lem}
\begin{proof}
    Let $k \in \N$; as in \cite{bernstein2018signsgd}, the strategy consists of breaking up the sum  in order to bound the both terms separately. 
    \begin{align*}
         \sum_{r = 0}^{k-1} (s_2^{i,r})^2  \prod_{t = r}^{k-1} (1-s_2^{i,t+1})^2 &=  \sum_{r = 0}^{\lfloor k/2 \rfloor - 1} (s_2^{i,r})^2  \prod_{t = r}^{k-1} (1-s_2^{i,t+1})^2 +  \sum_{r = \lfloor k/2 \rfloor }^{k-1} (s_2^{i,r})^2  \prod_{t = r}^{k-1} (1-s_2^{i,t+1})^2 \\
         &\leq (1-s_2^{i,k})^{2\lfloor k/2 \rfloor} \sum_{r = 0}^{\lfloor k/2 \rfloor - 1} (s_2^{i,r})^2 + (s_2^{i,\lfloor k/2 \rfloor - 1})^2 \sum_{r = \lfloor k/2 \rfloor }^{k-1} (1-s_2^{i,k})^{2(k-r-1)} \\
         & \leq  (s_2^{i,0})^2 \lfloor k/2 \rfloor (1-s_2^{i,k})^{2\lfloor k/2 \rfloor}  + \frac{8(s_2^{i,0})^2 }{k^{2\alpha_2}} \sum_{r = 0}^{\lfloor k/2 \rfloor} (1-s_2^{i,k})^{2r} \\
         &\leq  (s_2^{i,0})^2 k (1-s_2^{i,k})^{2 \lfloor k/2 \rfloor}  + \frac{8 (s_2^{i,0})^2 }{k^{2\alpha_2}(1 - (1- s_2^{i,k})^2)}\\
         &\leq  (s_2^{i,0})^2 k (1-s_2^{i,k})^{2\lfloor k/2 \rfloor}  + \frac{8 s_2^{i,0} }{k^{\alpha_2}(2 - s_2^{i,k})}.
    \end{align*}
Now, we are looking for $k$ such that
\begin{equation*}
    s_2^{i,0} k (1-s_2^{i,k})^{2\lfloor k/2 \rfloor} \leq \frac{1}{k^{\alpha_2}} \Leftrightarrow e^{2 \lfloor k/2 \rfloor \ln(1-s_2^{i,k})} \leq \frac{1}{(s_2^{i,0}) k^{1+\alpha_2}}.
\end{equation*}
As, $\ln(1-x) \leq -x$, it is sufficient to find $k$ such that 
\begin{align*}
    & \quad e^{-s_2^{i,0} \frac{k}{(k+1)^{\alpha_2}}} \leq \frac{1}{(s_2^{i,0}) k^{1+\alpha_2}} \\
    &\Leftrightarrow \frac{k}{(k+1)^{\alpha_2}} \geq \frac{\ln(s_2^{i,0}) + (1+\alpha_2) \ln(k)}{s_2^{i,0}}.
\end{align*}
Taking such a $k$ allows to complete the proof.
\end{proof}
Combining the three previous Lemmas allows bounding the variance term in the Proposition \ref{prop1}.
\begin{Prop}
    In the setting of Lemmas \ref{lem_var} and \ref{Lem_utile} and under Assumption \ref{assum1}.b, the variance term of Proposition \ref{prop1} is bounded by 
    \begin{equation*}
        \mathbb{E}[||\m^{i,k+1} - \Bar{\m}^{i,k+1}||_2^2] \leq  \frac{9 s_2^{i,0} L_0(F)^2 (n+4)^2 }{k^{\alpha_2}}+ o\left(\frac{1}{k^{\alpha_2}} \right).
    \end{equation*}  
    \label{prop2}
\end{Prop}
\begin{proof}
    By Lemmas \ref{Lem_decomp} and \ref{lem_var}, it follows that 
    \begin{align*}
        \mathbb{E}[||\m^{i,k+1} - \Bar{\m}^{i,k+1}||^2]  &\leq (s_2^{i,k})^2 \mathbb{E}[||\Tilde{\g}^{i,k} ||^2] + \sum_{r = 0}^{k-1} (s_2^{i,r})^2  \prod_{t = r}^{k-1} (1-s_2^{i,t+1})^2 \mathbb{E}[||\Tilde{\g}^{i,r} ||^2]\\
        &+ \prod_{t = 0}^{k} (1-s_2^{i,t})^2 \mathbb{E}[||\Tilde{\g}^{i,0} ||^2] \\
        &\leq \left ( (s_2^{i,k})^2  + \sum_{r = 0}^{k-1} (s_2^{i,r})^2  \prod_{t = r}^{k-1} (1-s_2^{i,t+1})^2 + \prod_{t = 0}^{k} (1-s_2^{i,t})^2  \right)  L_0(F)^2 (n+4)^2 .
    \end{align*}
    Now as $(s_2^{i,k})^2 = o\left( \frac{1}{k^{\alpha_2}} \right)$ and $\prod_{t = 0}^{k} (1-s_2^{i,t})^2 = o\left( \frac{1}{k^{\alpha_2}} \right) $, the result follows from Lemma \ref{Lem_utile}.
\end{proof}
\subsubsection{Bound on the bias term}
First, the bias term is bounded by a sum depending on $s_1^k$ and $s_2^k$. 
\begin{Lem}
    For the subproblem $i \in \N$ and at iteration $k \in \N$ of the algorithm \ref{algo1}, we have 
    \begin{equation*}
        \mathbb{E}[||\Bar{\m}^{i,k+1} - \nabla f^{\beta^i}(\x^{i,k}) ||_1] \leq 2n L_1(f^{\beta^i}) \left( \sum_{l = 0}^{k-1} s_1^{i,l} \prod_{t = l}^{k-1} (1-s_2^{i,t+1}) \right).
    \end{equation*}
    \label{Lem_reformulate}
\end{Lem}
\begin{proof}
    Foremost, observe that the quantity
    \begin{equation}
    \label{equa2}
        S^{i,k} :=  \left\{
    \begin{array}{ll}
        1 & \mbox{if } k = 0 \\
    s_2^{i,k}  + \sum_{r = 0}^{k-1} s_2^{i,r} \prod_{t = r}^{k-1} (1-s_2^{i,t+1})  + \prod_{t = 0}^k (1-s_2^{i,t})
          & \mbox{otherwise,}
    \end{array}
    \right.
    \end{equation}
    may be written recursively as 
    \begin{equation*}
        S^{i,k} = \left\{
    \begin{array}{ll}
        1 & \mbox{if } k = 0 \\
         s_2^{i,k} + (1-s_2^{i,k}) S^{i,k-1} & \mbox{otherwise.}
    \end{array}
\right.
    \end{equation*}
     Note that in its second expression $S^{i,k} = 1$ for all $k$. Therefore, by definition of $\Bar{m}_j^{i,k}$ and the previous result on $S^{i,k}$, it follows that
    \begin{align*}
        &\Bar{\m}^{i,k} = s_2^{i,k} \nabla f^{\beta^i}(\x^{i,k}) + \sum_{r = 0}^{k-1} s_2^{i,r} \prod_{t = r}^{k-1} (1-s_2^{i,t+1}) \nabla f^{\beta^i}(\x^{i,r}) + \prod_{t = 0}^k (1-s_2^{i,t}) \nabla f^{\beta^i}(\x^{i,0}) \\
        &\nabla f^{\beta^i}(\x^{i,k}) = \left( s_2^{i,k}  + \sum_{r = 0}^{k-1} s_2^{i,r} \prod_{t = r}^{k-1} (1-s_2^{i,t+1})  + \prod_{t = 0}^k (1-s_2^{i,t}) \right) \nabla f^{\beta^i}(\x^{i,k}).
    \end{align*}
   Thus
    \begin{align}
    \begin{split}
        \mathbb{E}[||\Bar{\m}^{i,k+1} - \nabla f^{\beta^i}(\x^{i,k}) ||_1] &\leq   \sum_{r = 0}^{k-1} s_2^{i,r} \prod_{t = r}^{k-1} (1-s_2^{i,t+1}) \mathbb{E}[||\nabla f^{\beta^i}(\x^{i,r}) - \nabla f^{\beta^i}(\x^{i,k})||_1] \\
        &+ \prod_{t = 0}^k (1-s_2^{i,t})  \mathbb{E}[||\nabla f^{\beta^i}(\x^{i,0}) - \nabla f^{\beta^i}(\x^{i,k})||_1].
    \end{split}
    \label{28}
    \end{align}
    By the smoothness of the function $f^{\beta^i}$, Lemma F.3 of \cite{bernstein2018signsgd} ensures that $\forall r \in [0, k-1]$
    \begin{align*}
        ||\nabla f^{\beta^i} (\x^{i,r}) - \nabla f^{\beta^i}(\x^{i, k})||_1 &\leq \sum_{l = r}^{k-1} ||\nabla f^{\beta^i}(\x^{i, l+1}) - \nabla f^{\beta^i}(\x^{i,l}) ||_1 \  \leq \ 2 n L_1(f^{\beta^i}) \sum_{l = r}^{k-1} s_1^{i,l}.
    \end{align*}
    Substituting this inequality in Equation (\ref{28}) gives 
 \begin{align}
    \begin{split}
        \mathbb{E}[||\Bar{m}^{i,k+1} - \nabla f^{\beta^i}(\x^{i,k}) ||_1] &\leq   2n L_1(f^{\beta^i}) S_1^{i,k} 
    \end{split}
    \end{align}
    where 
    \begin{equation*}
        S_1^{i,k} = \sum_{r = 0}^{k-1} s_2^{i,r} \sum_{l = r}^{k-1} s_1^{i,l}  \prod_{t = r}^{k-1} (1-s_2^{i,t+1}) 
        +  \sum_{l = 0}^{k-1} s_1^{i,l}  \prod_{t = 0}^k (1-s_2^{i,t}).
    \end{equation*}
    Reordering the terms in $S_1^k$, we obtain
    \begin{align*}
        S_1^{i,k} &= \sum_{l = 0}^{k-1} s_1^{i,l} \left( \sum_{r = 0}^l s_2^{i,r} \prod_{t = r}^{k-1} (1-s_2^{i,t+1}) + \prod_{t= 0}^k(1-s_2^{i,t}) \right) \\
        &= \sum_{l = 0}^{k-1} s_1^{i,l} \left( s_2^{i,l}\prod_{t = l}^{k-1} (1-s_2^{i,t+1}) +  \sum_{r = 0}^{l-1} s_2^{i,r} \prod_{t = r}^{k-1} (1-s_2^{i,t+1}) + \prod_{t = 0}^k(1-s_2^{i,t}) \right) \\
        &= \sum_{l = 0}^{k-1} s_1^{i,l} \prod_{t = l}^{k-1} (1-s_2^{i,t+1}) \underbrace{\left( s_2^{i,l} +  \sum_{r = 0}^{l-1} s_2^{i,r} \prod_{t = r}^{l-1} (1-s_2^{i,t+1}) + \prod_{t = 0}^{l}(1-s_2^{i,t}) \right)}_{S^{i,l} = 1} \\
        &= \sum_{l = 0}^{k-1} s_1^{i,l} \prod_{t = l}^{k-1} (1-s_2^{i,t+1}),
    \end{align*}
    which completes the proof.
\end{proof}
Second, the sum may be bounded by a term decreasing with $k$.
    \begin{Lem}
    For the subproblem $i \in \N$ and let $s_2^{i,k}  = \frac{s_2^{i,0}}{(k+1)^{\alpha_2}}$ and $s_1^{i,k} = \frac{s_1^{i,0}}{(k+1)^{\alpha_1}}$  with $s_1^{i,0} \in (0, 1), s_2^{i,0} \in (0,1)$ and $0 < \alpha_2 < \alpha_1 < 1$, then for $k$ such that 
    \begin{equation}
         \frac{k}{(k+1)^{\alpha_2}} \geq \frac{2 \left(\ln(s_2^{i,0}) + (1+\alpha_1 - \alpha_2)\ln(k)\right)}{s_2^{i,0}}
         \label{cond2}
    \end{equation}
    the following inequality holds
    \begin{equation}
        \sum_{l = 0}^{k-1} s_1^{i,l}  \prod_{t = l}^{k-1} (1-s_2^{i,t+1}) \leq \frac{5 s_1^{i,0} }{s_2^{i,0} k^{\alpha_1-\alpha_2}}.
    \end{equation}
    \label{Lem_utile_bis}
\end{Lem}

\begin{proof}
    The proof follows the proof of Lemma \ref{Lem_utile}. The sum is partitioned as follows:
    \begin{align*}
        \sum_{l = 0}^{k-1} s_1^{i,l} \prod_{t = l}^{k-1} (1-s_2^{i,t+1})  &= \sum_{l = 0}^{\lfloor k/2 \rfloor-1} s_1^{i,l} \prod_{t = l}^{k-1} (1-s_2^{i,t+1}) + \sum_{l = \lfloor k/2 \rfloor-1}^{k-1} s_1^{i,l} \prod_{t = l}^{k-1} (1-s_2^{i,t+1})\\
        &\leq (1-s_2^{i,k})^{\lfloor k/2 \rfloor} \sum_{l = 0}^{\lfloor k/2 \rfloor-1} s_1^{i,l} + s_1^{i, \lfloor k/2 \rfloor -1} \sum_{l = \lfloor k/2 \rfloor-1}^{k-1} (1-s_2^{i,k})^{k-r-1} \\
        &\leq s_1^{i,0} k (1-s_2^{i,k})^{\lfloor k/2 \rfloor}   + \frac{4 s_1^{i,0}}{k^{\alpha_1}(1-(1-s_2^{i,k}))} \\
        &=  \frac{s_1^{i,0} s_2^{i,0} k (1-s_2^{i,k})^{\lfloor k/2 \rfloor}}{ s_2^{i,0}}   + \frac{4 s_1^{i,0}}{s_2^{i,0} k^{\alpha_1 - \alpha_2}}.
    \end{align*}
    Now, as in Lemma \ref{Lem_utile} taking $k$ such that 
    \begin{equation*}
        \frac{k}{(k+1)^{\alpha_2}} \geq \frac{2 \left(\ln(s_2^{i,0}) + (1+\alpha_1 - \alpha_2)\ln(k)\right)}{s_2^{i,0}}
    \end{equation*}
    ensures that $s_2^{i,0} k (1-s_2^{i,k})^{\lfloor k/2 \rfloor} \leq \frac{1}{k^{\alpha_1 - \alpha_2}}$, which completes the proof.
\end{proof}
Finally, using the two previous Lemmas allows bounding the bias term.
\begin{Prop}
    In the setting of Lemma \ref{Lem_utile_bis}, the bias term of Proposition \ref{prop1} is bounded by 
    \begin{equation*}
        \mathbb{E}[||\Bar{\m}^{i,k+1} - \nabla f^{\beta^i}(\x^{i,k}) ||_1] \leq 10 n L_1(f^{\beta^i}) \frac{s_1^{i,0}}{s_2^{i,0} k^{\alpha_1-\alpha_2}}.
    \end{equation*}
    \label{prop3}
\end{Prop}
\begin{proof}
    The proof is a straightforward consequence of Lemmas \ref{Lem_reformulate} and \ref{Lem_utile_bis}.
\end{proof}
\subsubsection{Convergence in expectation of the ZOS algorithm}
As the different terms in the inequality of the Proposition \ref{prop1} have been bounded, the main result of this section may be derived in the following theorem.
\begin{Th}
For a subproblem $i \in \N$ and under Assumption \ref{assum1}, let $\alpha_1 \in (0, 1)$, $\alpha_2 \in (0, \alpha_1)$, $0<~s_1^{i,0}, s_2^{i,0} < 1$ and $K > C$ where $C \in \N$ satisfies Equations (\ref{cond1}) and (\ref{cond2}), we have 
\begin{align}
\begin{split}
    \mathbb{E}[||\nabla f^{\beta^i}(\x^{i,R})||_1] \leq \frac{1}{K^{1-\alpha_1} - \frac{C}{K^{\alpha_1}}} \Bigg( \frac{D_f^i}{s_1^{i,0}}  + \frac{n \sqrt{n} L_0(F) s_1^{i,0}}{\beta^i} \sum_{k = C}^K \frac{1}{k^{2 \alpha_1}} &+ 6 \sqrt{s_2^{i,0}} L_0(F) \sqrt{n}(n+4)  \sum_{k = C}^K \frac{1}{k^{\alpha_1 + \frac{\alpha_2}{2}}} \\
    & +  \frac{40 L_0(F) s_1^{i, 0} n \sqrt{n} }{s_2^{i,0} \beta^i} \sum_{k = C}^K \frac{1}{k^{2\alpha_1 - \alpha_2}} \Bigg),
\end{split}
\end{align}
where $f^{\beta^i}(\x^{i,C}) - \min_\x f^{\beta^i}(\x) \leq D_f^i$, $L_0(F)$ is the Lipschitz constant of $F$  and $R$ is randomly picked from a uniform distribution in $[C,K]$.
\label{conv_expect}
\end{Th}
\begin{proof}
    Let $C \in \N $ satisfying Equations (\ref{cond1}) and (\ref{cond2}) and sum over the inequality in Proposition \ref{prop1}, it follows that
    \begin{align*}
        \sum_{k = C}^K s_1^{i,k}  \mathbb{E}[||\nabla f^{\beta^i}(\x^{i,k})||_1] &\leq \mathbb{E}[f^{\beta^i}(x^{i,C}) - f^{\beta^i}(\x^{i,K+1})] + \frac{n L_1(f^{\beta^i})}{2} \sum_{k = C}^K (s_1^{i,k})^2 \\
        &+ 2 \sqrt{n} \sum_{k = C}^K s_1^{i,k}  \sqrt{\mathbb{E}[||\m^{i,k+1} - \Bar{\m}^{i,k+1}||_2^2]}\\
        &+ 2 \sum_{k = C}^K s_1^{i,k} \mathbb{E}[||\Bar{\m}^{i,k+1} - \nabla f^{\beta^i}(\x^{i,k}) ||_1].
    \end{align*}
    By substituting the results of Proposition \ref{prop2} and \ref{prop3} in the previous inequality, we obtain
    \begin{align*}
       \sum_{k = C}^K s_1^{i,k}  \mathbb{E}[||\nabla f^{\beta^i}(\x^{i,k})||_1] &\leq \mathbb{E}[f^{\beta^i}(x^{i,C}) - f^{\beta^i}(\x^{i,K+1})] + \frac{nL_1(f^{\beta^i})}{2} \sum_{k = C}^K (s_1^{i,k})^2 \\
       &+ 6\sqrt{s_2^{i,0}} L_0(F) (n+4) \sqrt{n}  \sum_{k = C}^K \frac{s_1^{i,0}}{k^{\alpha_1 + \frac{\alpha_2}{2}}} 
    +  \frac{20  L_1(f^{\beta^i})  s_1^{i, 0} n }{s_2^{i,0}} \sum_{k = C}^K \frac{s_1^{i,0}}{k^{2\alpha_1 - \alpha_2}}.
    \end{align*}
    Dividing both sides by $s_1^{i,0} K^{-\alpha_1} (K - C)$, picking $R$ randomly uniformly in $[C, K]$ and using the definition of $D_f^i$ given that $\min_\x f(\x) \leq f(\x)$ for all $\x$, we get
    \begin{align*}
        \mathbb{E}[||\nabla f^{\beta^i}(x^{i,R})||_1] &= \frac{1}{K-C} \sum_{k = C}^K \mathbb{E}[||\nabla f^{\beta^i}(\x^{i,k})||_1]
        \leq \frac{1}{K-C} \sum_{k = C}^K \frac{K^{\alpha_1}}{k^{\alpha_1}} \mathbb{E}[||\nabla f^{\beta^i}(\x^{i,k})||_1]\\
        &\leq \frac{1}{K^{1-\alpha_1} - \frac{C}{K^{\alpha_1}}} \Bigg( \frac{D_f^i}{s_1^{i,0}}  + \frac{n L_1(f^{\beta^i}) s_1^{i,0}}{2} \sum_{k = C}^K \frac{1}{k^{2 \alpha_1}} + 6 \sqrt{s_2^{i,0}} L_0(F) (n+4) \sqrt{n}  \sum_{k = C}^K \frac{1}{k^{\alpha_1 + \frac{\alpha_2}{2}}} \\
    & +  \frac{20 L_1(f^{\beta^i}) s_1^{i, 0} n }{s_2^{i,0}} \sum_{k = C}^K \frac{1}{k^{2\alpha_1 - \alpha_2}} \Bigg).
    \end{align*} 
Recalling that $L_1(f^{\beta^i}) = \frac{2 \sqrt{n} L_0(F)}{\beta^i}$ (see \cite[Lemma 2]{nesterov2017random}) completes the proof.
\end{proof}
This theorem allows proving the convergence in expectation of the norm of the gradient when $\alpha_1$ and $\alpha_2$ are chosen adequately. In particular, the following corollary provides the convergence when $\alpha_1 = \frac{3}{4}$ and $\alpha_2 = \frac{1}{2}$. 
\begin{Cor}
\label{cor1}
    Under the same setting of Theorem with $\beta^i \approx 1$  $\alpha_1 = \frac{3}{4}$, $\alpha_2 = \frac{1}{2}$, $s_1^{i,0} = \frac{1}{n^{\frac{3}{4}}}$ and $s_2^{i,0} \approx 1 $, we have 
    \begin{equation}
        \mathbb{E}[||\nabla f^{\beta^i}(\x^{i,R}) ||_2] = O \left( \frac{n^{\frac{3}{2}}}{K^{1/4}} \ln(K)\right).
        \label{equa_conv_expect}
    \end{equation}
\end{Cor}
\begin{proof}
    The result is a direct consequence of Theorem \ref{conv_expect} with the specified constant and by noting that $||\cdot||_2 \leq ||\cdot||_1$ in $\R^n$.
\end{proof}
In \cite{chen2019zo, ghadimi2013stochastic, liu2018signsgd}, the function $F$ is assumed to be smooth with Lipschitz continuous gradient. 
In the present work, $F$ is only assumed to be Lipschitz continuous. 
This has two main consequences on the result of convergence: the dependence of the dimension on the convergence rate is larger. Furthermore, while $\beta$ must be chosen relatively small in the smooth case, it is interesting to note that it does not have to be this way in the nonsmooth case.

\subsubsection{The convex case}

The convergence results of the ZOS algorithm has been derived in the non-convex case. In the next theorem, convergence results are derived when the function $f^{\beta^i}$ is convex. 
\begin{Th}
\label{th_convex}
Under Assumption \ref{assum1}, suppose moreover that $f^{\beta^i}$ is convex and there exists $\rho$ such that $\rho = \max_{k \in \N} ||\x^{i,k} - \x^{i, *}||$, then by setting 
\begin{equation}
    s_1^{i,k} = \frac{2 \rho}{(k+1)},\; s_2^{i,k} = \frac{1}{(k+1)^{\frac{2}{3}}} \text{ and } \Gamma^k := \prod_{l = 2}^k \left(1 - \frac{2}{k+1} \right) = \frac{2}{k(k+1)} \text{ with } \Gamma^1 = 1,
    \label{setting}
\end{equation}
it follows that 
\begin{equation}
    \mathbb{E}[f^{\beta^i}(\x^{i, K}) - f^{\beta^i}(\x^{*})] \leq \frac{4\rho^2 n \sqrt{n}L_0(F)}{\beta^i K^{\frac{1}{3}}}.
    \label{res_conv_f}
\end{equation}
and 
\begin{equation}
    \mathbb{E}[||\nabla f^{\beta^i}(\x^{i,R})||] \leq \frac{2L_0(F)}{K^2} + \frac{4 \rho n\sqrt{n}L_0(F)}{\beta^i K^{\frac{1}{3}}}
    \label{res_conv_g}
\end{equation}
where $R$ is a random variable in $[0,K-1]$ whose the probability distribution is given by
\begin{equation*}
    \mathbb{P}(R~=~k)~=~\frac{s_1^{i,k}/\Gamma^{k+1}}{\sum_{k = 0}^{K-1} s_1^{i,k}/\Gamma^{k+1}}.
\end{equation*}
\end{Th}
\begin{proof}
Under the assumptions in the statement of the Theorem, it follows by Proposition \ref{prop1} that 
\begin{align}
\begin{split}
    \mathbb{E}[f^{\beta^i}(\x^{i, k+1}) - f^{\beta^i}(\x^{i, *})] &\leq \mathbb{E}[f^{\beta^i}(\x^{i, k}) - f^{\beta^i}(\x^{i, *})] -s_1^{i,k}\mathbb{E}[||\nabla f^{\beta^i}(\x^{i,k})||] + \frac{nL_1(f^{\beta^i})}{2} (s_1^{i,k})^2 \\
    &+ 2 s_1^{i,k} \mathbb{E}[||\Bar{\m}^{i,k+1}  - \nabla f^{\beta^i}(\x^{i,k})||_1] + 2s_1^{i,k}  \sqrt{n} \sqrt{\mathbb{E}[||\m^{i,k+1} - \Bar{\m}^{i,k+1}||^2]} \\
      &\leq \mathbb{E}[f^{\beta^i}(\x^{i, k}) - f^{\beta^i}(\x^{i, *})] -s_1^{k}\mathbb{E}[||\nabla f^{\beta^i}(\x^{i,k})||] + \frac{ 4 \rho^2 n\sqrt{n}L_0(F) }{\beta^i (k+1)^{\frac{4}{3}}},
      \end{split}
      \label{equa10}
\end{align}
where the last inequality follows thanks to Propositions \ref{prop2}, \ref{prop3} with $L_1(f^{\beta^i})= \frac{2L_0(F) \sqrt{n}}{\beta^i} $ and the values of  $s_1^{i,k} $ and $s_2^{i,k}$. Now, by convexity assumption of $f^{\beta^i}$ and the bound $\rho$, the following holds
\begin{align*}
    f^{\beta^i}(\x^{i,k}) - f^{\beta^i}(\x^{i,*}) &\leq \nabla f^{\beta^i}(\x^{i,k})^T (\x^{i,k} - \x^{i,*}) \\
    &\leq ||\nabla f^{\beta^i}(\x^{i,k})|| \, ||\x^{i,k} - \x^{i,*}||\\
    &\leq \rho ||\nabla f^{\beta^i}(\x^{i,k})||.
\end{align*}
Thus, by substituting this result into Equation (\ref{equa10}), it follows that 
\begin{equation*}
    \mathbb{E}[f^{\beta^i}(\x^{i, k+1}) - f^{\beta^i}(\x^{i, *})] \leq \left( 1 - \frac{2}{(k+1)} \right)\mathbb{E}[f^{\beta^i}(\x^{i, k}) - f^{\beta^i}(\x^{i, *})]  + \frac{ 4\rho^2 n \sqrt{n}L_0(F)} {\beta^i (k+1)^{\frac{4}{3}}}.
\end{equation*}
Now by dividing by $\Gamma^{k+1}$ both sides of the equation and summing up the inequalities, it follows that
\begin{align*}
    \frac{\mathbb{E}[f^{\beta^i}(\x^{i, K}) - f^{\beta^i}(\x^{i, *})]}{\Gamma^{K}} &\leq \frac{ 4\rho^2 n \sqrt{n}L_0(F)}{\beta^i} \sum_{k = 0}^{K-1} \frac{1}{\Gamma^{k+1} (k+1)^{\frac{4}{3}}}\\
    &\leq \frac{ 4\rho^2 n \sqrt{n}L_0(F)}{\beta^i} \sum_{k = 0}^{K-1} (k+1)^{\frac{2}{3}}. 
\end{align*}
Thus
\begin{equation*}
    \mathbb{E}[f^{\beta^i}(\x^{i, K}) - f^{\beta^i}(\x^{i, *})] \leq  \frac{ 4\rho^2 n \sqrt{n}L_0(F)}{\beta^i} \Gamma^K \sum_{k = 0}^{K-1} (k+1)^{\frac{2}{3}} \leq  \frac{ 4\rho^2 n \sqrt{n}L_0(F)}{\beta^i K^{\frac{1}{3}}} .
\end{equation*}

Now, the second part of the proof may be demonstrated. By Equation (\ref{equa10}), it follows also that
\begin{equation*}
    s_1^{i,k}\mathbb{E}[||\nabla f^{\beta^i}(\x^{i,k})||] \leq  \mathbb{E}[f^{\beta^i}(\x^{i, k}) - f^{\beta^i}(\x^{i, *})] - \mathbb{E}[f^{\beta^i}(\x^{i, k+1}) - f^{\beta^i}(\x^{i, *})] + \frac{  4\rho^2 n\sqrt{n}L_0(F) }{\beta^i (k+1)^{\frac{4}{3}}}.
\end{equation*}

As in the previous part, by dividing both sides by $\Gamma^{k+1}$, summing up the inequalities and noting $\Bar{f^k} = \mathbb{E}[f^{\beta^i}(\x^{i, k}) - f^{\beta^i}(\x^{i, *})]$, we obtain
\begin{equation*}
    \sum_{k = 0}^{K-1} \frac{s_1^{i,k} }{ \Gamma^{k+1}} \mathbb{E}[||\nabla f^{\beta^i}(\x^{i,k})||] \leq \sum_{k = 0}^{K-1} \frac{\Bar{f^k} -\Bar{f^{k+1}}}{\Gamma^{k+1}} +  \frac{ 4\rho^2 n \sqrt{n}L_0(F)}{\beta^i} \sum_{k = 0}^{K-1}\frac{1}{\Gamma^{k+1} (k+1)^{\frac{4}{3}}}.
\end{equation*}
Then, again by dividing both sides by $\sum_{k = 0}^{K-1} \frac{s_1^{i,k}}{\Gamma^{k+1}}$ it follows that 
\begin{align*}
    \mathbb{E}[||\nabla f^{\beta^i}(\x^{i,R})||] &= \frac{\sum_{k = 0}^{K-1} \frac{s_1^{i,k} }{ \Gamma^{k+1}} \mathbb{E}[||\nabla f^{\beta^i}(\x^{i,k})||]}{\sum_{k = 0}^{K-1} \frac{s_1^{i,k}}{\Gamma^{k+1}}} \\
    &\leq \frac{1}{\sum_{k = 0}^{K-1} \frac{s_1^{i,k}}{\Gamma^{k+1}}} \left( \sum_{k = 0}^{K-1} \frac{\mathbb{E}[\Bar{f^k} -\Bar{f^{k+1}}]}{\Gamma^{k+1}} + \frac{ 4\rho^2 n \sqrt{n}L_0(F)}{\beta^i} \sum_{k = 0}^{K-1}\frac{1}{\Gamma^{k+1} (k+1)^{\frac{4}{3}}} \right),
\end{align*}
where $R$ is a random variable whose the distribution is given in the statement of the theorem. Now, as in Equation (2.21) of \cite{balasubramanian2022zeroth}, the following inequalities hold
\begin{eqnarray*}
        \sum_{k = 0}^{K-1} \frac{\Bar{f^k} -\Bar{f^{k+1}}}{\Gamma^{k+1}} \leq \Bar{f^0} + \sum_{k =  1}^{K-1} \frac{2}{\Gamma^{k+1}(k+1)} \Bar{f^k}
    & \quad \mbox{ and } \quad&
    \sum_{k = 0}^{K-1} \frac{s_1^{i,k}}{\Gamma^{k+1}} = \frac{\rho}{\Gamma^K}.
\end{eqnarray*}
Thus, by substituting these in the inequality involving the expectation, we obtain
\begin{align*}
    \mathbb{E}[||\nabla f^{\beta^i}(\x^{i,R})||] &\leq \frac{\Gamma^K}{\rho}\left( \mathbb{E}[\Bar{f^0}] + \sum_{k =  1}^{K-1} \frac{2}{\Gamma^{k+1}(k+1)} \mathbb{E}[\Bar{f^k}] + \frac{ 4\rho^2 n \sqrt{n}L_0(F)}{\beta^i} \sum_{k = 0}^{K-1}\frac{1}{\Gamma^{k+1} (k+1)^{\frac{4}{3}}}\right)\\
    &\leq \frac{\Gamma^K}{\rho}\left( \mathbb{E}[\Bar{f^0}]  + \frac{8 \rho n\sqrt{n}L_0(F)}{\beta^i} \sum_{k = 0}^{K-1}\frac{1}{\Gamma^{k+1} (k+1)^{\frac{4}{3}}}\right)\\
    &\leq \frac{2 L_0(F)}{K^2} + \frac{8 \rho n\sqrt{n}L_0(F)}{\beta^i K^{\frac{1}{3}}},
\end{align*}
where the second inequality follows from Equation (\ref{res_conv_f}).
\end{proof}

\subsubsection{Summary of convergence rates and complexity guarantees}
The result obtained in Equation (\ref{equa_conv_expect}) is consistent with the convergence results of other ZO methods. To gain a better understanding of its performance, this result is compared  with those of four other algorithms from different perspectives: the assumptions, the measure used, the convergence rate, and the function query complexity. All  methods seek a solution to a  stochastic optimization problem; the comparison is presented in Table \ref{tab1}. Since the convergence rate of the ZO-Signum and ZO-signSGD algorithms is measured using $||\nabla f(\x)||$, but $||\nabla f(x)||^2$ is used by ZO-adaMM and ZO-SGD, Jensen's inequality is used to rewrite  convergence rates in terms of gradient norm. 
\begin{itemize}
    \item for ZO-SGD \cite{ghadimi2013stochastic}
    \begin{equation*}
        \mathbb{E}[||\nabla f(\x)||] \leq \sqrt{\mathbb{E}[||\nabla f(\x)||^2]} \leq \sqrt{O\left(\frac{\sigma \sqrt{n}}{\sqrt{K}} + \frac{n}{K} \right)} \leq O\left( \frac{\sqrt{\sigma}n^{\frac{1}{4}}}{K^{\frac{1}{4}}} + \frac{\sqrt{n}}{\sqrt{K}}\right),
        \end{equation*}
    \item for ZO-adaMM \cite{chen2019zo}
    \begin{align*}
        \mathbb{E}[||\nabla f(\x)||] &\leq \sqrt{\mathbb{E}[||\nabla f(\x)||^2]} \leq \sqrt{O\left( \left(\frac{n}{\sqrt{K}} + \frac{n^2}{K}\right)\sqrt{\ln(K) + \ln(n)}  \right)}\\
        &\leq O \left( \left( \frac{\sqrt{n}}{K^{\frac{1}{4}}} + \frac{n}{\sqrt{K}}\right) \left( \ln(K) + \ln(n) \right)^{\frac{1}{4}} \right),
        \end{align*}
\end{itemize}
where the third inequalities are due to $\sqrt{a^2 + b^2} \leq a + b$, for $a,b \geq 0$. For ZO-signSGD, unless the value of $b$  depends on $K$, the algorithm's convergence is only guaranteed within some ball around the solution, making it difficult to compare with other methods. Thus, in the non-convex case, after this transformation, it becomes apparent that ZO-Signum has a convergence rate of $O\left(\frac{n^{\frac{3}{4}}}{\sqrt{\sigma}} \right)$ and $O(\sqrt{n})$ worse than that of ZO-SGD and ZO-adaMM, respectively. This may be attributed to the milder assumption made on the function $F$ in the present work, which also explains why the convergence is relative to $f^{\beta}$. In the convex case, ZO-Signum has a convergence rate of $O\left(\frac{n K^{\frac{1}{6}}}{\sigma} \right)$
worse than ZSCG and $O\left(\sqrt{n} K^{\frac{1}{6}} \right)$ worse than ZO-SGD. This may be explained because the $sign(\cdot)$ operator looses the magnitude information of the gradient when it applied. 
This problem may be fixed as in \cite{karimireddy2019error} but it outside the scope  of this work. Finally, all methods but ZO-signSGD are momentum-based versions of the original ZO-SGD method. Although the momentum-based versions are mostly used in practice, it is interesting to notice that none of these methods possess a better convergence rate than the original ZO-SGD method. The next section provides some clues on the interests of the momentum-based method.

\newpage
\begin{table}[htb!]
    \centering
    \footnotesize
        \caption{Summary of convergence rate and query complexity of various ZO-algorithms given K iterations.}
\begin{tabular}{|c|c|c|c|c|}
        \hline
          Method & Assumptions & Measure & Convergence rate & Queries\\
         \hline
         \multirow{2}{*}{ZO-SGD \cite{ghadimi2013stochastic}} & \multirow{1}{*}{$F(\cdot, \boldsymbol{\xi}) \in \mathcal{C}^{1 +}$} & \multirow{2}{*}{$\mathbb{E}[||\nabla f(\x^R)||_2]$} &\multirow{2}{*}{$O\left( \frac{\sqrt{\sigma}n^{\frac{1}{4}}}{K^{\frac{1}{4}}} + \frac{\sqrt{n}}{\sqrt{K}}\right)$}  & \multirow{2}{*}{$O(K)$}\\
         &$\mathbb{E}[||\nabla F(\x, \boldsymbol{\xi}) - \nabla f(\x)||^2] \leq \sigma^2$ & & &\\
         \hline
         \multirow{2}{*}{ZO-signSGD \cite{liu2018signsgd}} & \multirow{1}{*}{$F(\cdot, \boldsymbol{\xi}) \in \mathcal{C}^{0 +}$} & \multirow{3}{*}{$\mathbb{E}[||\nabla f(\x^R)||_2]$} &\multirow{3}{*}{$O\left(\frac{\sqrt{n}}{\sqrt{K}} + \frac{\sqrt{n}}{\sqrt{b}} + \frac{n}{\sqrt{bq}} \right)$}  & \multirow{3}{*}{$O(bqK)$}\\
         &$F(\cdot, \boldsymbol{\xi}) \in \mathcal{C}^{1 +}$ & & &\\
         &$||\nabla F(\x, \boldsymbol{\xi})||_2 \leq \eta$& & & \\
         \hline
         \multirow{2}{*}{ZO-adaMM \cite{chen2019zo}} & \multirow{1}{*}{$F(\cdot, \boldsymbol{\xi}) \in \mathcal{C}^{0 +}$} & \multirow{3}{*}{$\mathbb{E}[||\nabla f(\x^R)||_2]$} &\multirow{3}{*}{$O \left( \left( \frac{\sqrt{n}}{K^{\frac{1}{4}}} + \frac{n}{\sqrt{K}}\right) \left( \ln(K) + \ln(n) \right)^{\frac{1}{4}} \right) $}  & \multirow{3}{*}{$O(K)$}\\
         &  \multirow{1}{*}{$F(\cdot, \boldsymbol{\xi}) \in \mathcal{C}^{1 +}$} & & &\\
         &$||\nabla F(\x, \boldsymbol{\xi})||_\infty \leq \eta$& & & \\
         \hline
         \textbf{ZO-Signum} & $F(\cdot, \boldsymbol{\xi}) \in \mathcal{C}^{0 +}$ & $\mathbb{E}[||\nabla f^{\beta}(\x^R)||_2]$ & $O\left(\frac{n \sqrt{n}}{K^{\frac{1}{4}}} \ln(K) \right)$ & $O(K)$\\
         \hline
         \textbf{ZO-Signum} & $F(\cdot, \boldsymbol{\xi}) \in \mathcal{C}^{0 +}$, 
         $f$ convex & $\mathbb{E}[f^{\beta^i}(\x^{i,K}) - f^{\beta^i}(\x^{i,*})]$ & $O\left(\frac{n\sqrt{n}}{K^{\frac{1}{3}}}  \right)$ & $O(K)$\\
          \hline
          ZO-SGD \cite{nesterov2017random} & $F(\cdot, \boldsymbol{\xi}) \in \mathcal{C}^{0 +}$, $f$ convex & $\mathbb{E}[f(\x^{i,K}) - f(\x^{i,*})]$ & $O\left(\frac{n}{\sqrt{K}} \right)$ & $O(K)$ \\
         \hline
          \multirow{2}{*}{Modified ZSCG \cite{balasubramanian2022zeroth}} & \multirow{1}{*}{$F(\cdot, \boldsymbol{\xi}) \in \mathcal{C}^{1 +}$, $F$ convex} & \multirow{2}{*}{$\mathbb{E}[f(\x^{i,K}) - f(\x^{i,*})]$} &\multirow{2}{*}{$O\left( \frac{\sigma \sqrt{n}}{\sqrt{K}} \right)$}  & \multirow{2}{*}{$O(K)$}\\
         &$\mathbb{E}[||\nabla F(\x, \boldsymbol{\xi}) - \nabla f(\x)||^2] \leq \sigma^2$ & & &\\
         \hline
    \end{tabular}
    \label{tab1}
\end{table}

\subsection{Convergence of the SSO algorithm}

The convergence analysis from the previous subsection is in expectation, i.e., it establishes the expected convergence performance over many executions of the ZO-Signum algorithm. As in \cite{ghadimi2013stochastic}, we now focus on the performance of a single run.  A second hierarchical workflow of the different theoretical results is presented in Table \ref{tab_conv_sso}.
\begin{table}[htb!]
    \centering
    \footnotesize
        \caption{Workflow of Lemmas/Propositions/Theorems for the SSO convergence analysis}
        \label{tab_conv_sso}
    \begin{tabular}{|C{2cm}|c|c|c|c|C{2cm}|C{2.1cm}|}
        \hline
          Assumptions on $F$ & \multicolumn{3}{c|}{Preliminary results} & Intermediate results & Main result & When $f^\beta$ is convex  \\
         \hline
         \multirow{2}{*}{Assumptions} & \multirow{1}{*}{Lemma \ref{lem_411}}& & &\multirow{3}{*}{Lemma \ref{bound_prob1} }& \multirow{5}{*}{Theorem \ref{main_th} (i)}  &  \multirow{5}{*}{Theorem \ref{main_th} (ii)}\\ 
         \cline{2-2}
          \multirow{2}{*}{\ref{assum1}, \ref{assum2} and \ref{assum3}} & \multirow{1}{*}{Proposition \ref{prop3}}  & \multirow{1}{*}{Lemma \ref{bound_grad}} &  \multirow{1}{*}{Lemma  \ref{lem_convex}} & &&\\
           \cline{2-2} 
          \multirow{2}{*}{which imply}& \multirow{1}{*}{Lemma \ref{lem_trunc_gaus}.3 }& & & & &\\
          \cline{2-4} 
         \multirow{2}{*}{$L_1(f^{\beta^i}) \leq L_1(f)$} &   \multicolumn{3}{c|}{Theorem \ref{conv_expect}} & & &\\
         \cline{2-5} 
         &   \multicolumn{3}{c|}{Proposition \ref{prop2}} &  \multirow{2}{*}{Lemma \ref{bound_prob2} } & &\\
         \cline{2-4} 
         &  \multicolumn{3}{c|}{Proposition \ref{prop3}} & & &\\
         \hline
    \end{tabular}
\end{table}
Unlike \cite{ghadimi2013stochastic},  our analysis is  based on a sequential optimization framework
 rather than a post-optimization process.
Our SSO algorithm uses the norm of the momentum as an indicator of the quality of the current solution. In order to analyze the rate of convergence of this algorithm, the following additional assumptions are made on the function $F$. The first assumption concerns the smoothness of the function $F$. The assumption of smoothness is used only to guarantee that $L_1(f^{\beta^i})$ is a constant with respect to $\beta^i$, contrarily to the non-smooth case (see \cite[Equation (12)]{nesterov2017random}).
\begin{assumption}
    The function $F(\cdot, \boldsymbol{\xi})$ has $L_1(F)$-Lipschitz continuous gradient.
    \label{assum2}
\end{assumption}
The second assumption concerns the local convexity of the function $f^{\beta}$. 
\begin{assumption}
    Let  $(\x^{i,0})$ be a sequence of points produced by Algorithm \ref{algo2} and $\x^{i,*}$ a sequence of local minima of $f^{\beta^i}$. We assume that there exists a threshold $I \in \N$ and a radius $\rho >0$ such that $\forall i \geq I$:
    \begin{enumerate}
        \item $f^{\beta^i}$ is convex on the ball $\mathcal{B}_{\rho}(\x^{i,*}) := \{ \x \in \R^n: ||x - x^{i,*}|| < \rho \}$;
        \item $\x^{i,0} \in \mathcal{B}_{\rho}(\x^{i,*})$.
    \end{enumerate}
    \label{assum3}
\end{assumption}
Under these assumptions, we will prove that if the norm of the momentum vector $\m$ is below some threshold, then this threshold can be used to bound the norm of the gradient. Second, an estimate for the number of iterations required to reduce the norm of $\m$ below the threshold is provided. The next lemma is simply technical and  demonstrates the link between $\Bar{\m}$ and $\m$.
\begin{Lem}
    For any subproblem $i \in \N$ and iteration $k \geq 1$, the following equality holds
    \begin{equation*}
        \mathbb{E}[\m^{i,k} | \x^{i,k-1}] = \mathbb{E}[\Bar{\m}^{i,k} | \x^{i,k-1}],
    \end{equation*}
    where $\Bar{\m}^{i,k}$ is defined recursively in Proposition \ref{prop1}.
    \label{lem_411}
\end{Lem}
\begin{proof}
    The proof is conducted by induction on $k$. For $k = 1$, setting $\m^{i,0} = \Tilde{\nabla} f^{\beta^i}(\x^{i,0}, \boldsymbol{\xi}^0)$ implies
    \begin{align*}
        \m^{i,1} &= s_2^{i,0} \Tilde{\nabla} f^{\beta^i}(\x^{i,0}, \boldsymbol{\xi}^0) + (1-s_2^{i,0}) \m^{i,0} = \Tilde{\nabla} f^{\beta^i}(\x^{i,0}, \boldsymbol{\xi}^0).
    \end{align*}
    In the same way, $\Bar{\m}^{i,1} = \nabla f^{\beta^i}(\x^{i,0})$. Therefore, we have
    \begin{equation*}
        \mathbb{E}[\m^{i,1} | \x^{i,0}] = \mathbb{E} [\Tilde{\nabla} f^{\beta^i}(\x^{i,0}, \boldsymbol{\xi}^0)| \x^{i,0}] = \nabla f^{\beta^i}(\x^{i,0}) = \mathbb{E}[\nabla f^{\beta^i}(\x^{i,0}) | \x^{i,0}] = \mathbb{E}[\Bar{\m}^{i,1} | \x^{i,0}].
    \end{equation*}
    Now, suppose that the induction assumption is true for a given $k \in \N$, then 
    \begin{align*}
        \mathbb{E}[\m^{i, k+1} | \x^{i,k}] = s_2^{i,k} \nabla f^{\beta^i}(\x^{i,k}) + (1-s_2^{i,k}) \mathbb{E}[\m^{i,k} | \x^{i,k}].
    \end{align*}
    Now, by the law of total expectation
    \begin{align*}
        \mathbb{E}[\m^{i,k} |\x^{i,k}] &= \mathbb{E}[ \mathbb{E}[\m^{i,k}|\x^{i,k}, \x^{i,k-1}] |\x^{i,k}]\\
        &= \mathbb{E}[ \mathbb{E}[\m^{i,k}| \x^{i,k-1}] |\x^{i,k}]\\
        &= \mathbb{E}[ \mathbb{E}[\Bar{\m}^{i,k}| \x^{i,k-1}] |\x^{i,k}] \quad  \text{(by the induction assumption)}\\
        &= \mathbb{E}[ \Bar{\m}^{i,k}|  \x^{i,k}].
    \end{align*}
    Thus as $\mathbb{E}[\nabla f^{\beta^i}(\x^{i,k}) |\x^{i,k}] = \nabla f^{\beta^i}(\x^{i,k})$, it follows that
    \begin{align*}
        \mathbb{E}[\m^{i, k+1} | \x^{i,k}] &= s_2^{i,k} \nabla f^{\beta^i}(\x^{i,k}) + (1-s_2^{i,k}) \mathbb{E}[\m^{i,k} | \x^{i,k}]\\
        &= s_2^{i,k} \mathbb{E}[\nabla f^{\beta^i}(\x^{i,k}) |\x^{i,k}]+ (1-s_2^{i,k}) \mathbb{E}[\Bar{\m}^{i,k} | \x^{i,k}]\\
        &= \mathbb{E}[\Bar{\m}^{i,k+1} | \x^{i,k}],
    \end{align*}
    which completes the proof.
\end{proof}
The following lemma shows that if $||\m||$ is below a certain threshold, then this threshold can be used to bound the norm of the gradient. 
\begin{Lem}
    For a subproblem $i \in \N$, let $K_i \in \N$ denote the first iteration in Algorithm \ref{algo1} for which $||\m^{i,K_i}|| \leq \frac{L\beta^i}{4 \beta^0}$ , then, under Assumption \ref{assum3} the norm of the gradient of the function $f^{\beta^i}$ at $\x^{i,K}$  may be bounded as follows
\begin{align*}
    || \nabla f^{\beta^i}(\x^{i,K_i}) || \leq \frac{L \beta^i}{4 \beta^0} + 10nL_1(F) \frac{s_1^{i,0}}{s_2^{i,0} K_i^{\alpha_1 - \alpha_2}}.
\end{align*}
Moreover, if the problem $i+1$ is considered, the gradient of the function $f^{\beta^{i+1}}$ may be bounded at the point $\x^{i,K_i} = \x^{i+1, 0}$ as follows
\begin{equation*}
 ||\nabla f^{\beta^{i+1}}(\x^{i+1, 0}) || \leq || \nabla f^{\beta^i}(\x^{i,K_i}) || + L_1(F) (n+3)^{\frac{3}{2}} (\beta^i - \beta^{i+1}).
\end{equation*}

\label{bound_grad}
\end{Lem}
\begin{proof}
   Let $K_i$ be taken as in the statement of the lemma. The norm of the gradient may be bounded as follows,
\begin{align*}
    ||\nabla f^{\beta^i}(\x^{i,K_i})|| &\leq ||\mathbb{E}[\m^{i,K_i} | \x^{i,K_i}] || + ||\nabla f^{\beta^i}(\x^{i,K_i}) - \mathbb{E}[\m^{i,K_i} | \x^{i,K_i}] ||\\
    &\leq \mathbb{E}[ ||\m^{i,K_i}|| \; | \x^{i,K_i}] + ||\nabla f^{\beta^i}(\x^{i,K_i}) - \mathbb{E}[ \Bar{\m}^{i,K_i} | \x^{i,K_i}] ||,  \\
\end{align*}
where the second inequality follows from Jensen's inequality and Lemma \ref{lem_411}. Now, using $||\m^{i,K_i}|| \leq \frac{L\beta^i}{4 \beta^0}$, $\mathbb{E}[\nabla f^{\beta^i}(\x^{i,K}) |\x^{i,K_i}] = \nabla f^{\beta^i}(\x^{i,K_i})$, $L_1(f^{\beta^i}) \leq L_1(F)$ and the result of Proposition \ref{prop3} completes the first part of the proof
\begin{align*}
    ||\nabla f^{\beta^i}(\x^{i,K_i})|| &\leq \frac{L \beta^i}{4 \beta^0} + \mathbb{E}[ ||\nabla f^{\beta^i}(\x^{i,K_i}) -  \Bar{\m}^{i,K_i}|| \; | \x^{i,K_i}] \\
    &\leq \frac{L \beta^i}{4 \beta^0} + 10nL_1(F) \frac{s_1^{i,0}}{s_2^{i,0} K_i^{\alpha_1 - \alpha_2}}.
\end{align*}
The second part of the proof follows directly by applying the triangular inequality and the result in Lemma \ref{lem_trunc_gaus}.3 because $\x^{i, K_i} = \x^{i+1, 0}$.
\end{proof}

Under Assumption \ref{assum2}, the expected difference between the values of $f^{\beta^i}$ at $\x^{i,0}$ and its optimal value is bounded in the next Lemma.
\begin{Lem}
    Let $I$ be the threshold from Assumption \ref{assum2}. If $i \geq I$, then
    \begin{equation}
        \mathbb{E}[f^{\beta^{i+1}}(\x^{i+1,0}) - f^{\beta^{i+1}}(\x^{i+1,*})] \leq \rho \left( \frac{L \beta^i}{4 \beta^0} + 10nL_1(F) \frac{s_1^{i,0}}{s_2^{i,0} K_i^{\alpha_1 - \alpha_2}} +  L_1(F) (n+3)^{\frac{3}{2}} (\beta^i - \beta^{i+1}) \right).
    \end{equation}
    \label{lem_convex}
\end{Lem}
\begin{proof}
    Convexity of the function $f^{\beta^i}$ on the ball $\mathcal{B}_{\rho}(\x^{i,*})$ implies
    \begin{align*}
        \mathbb{E}[f^{\beta^{i+1}}(\x^{i+1,0}) - f^{\beta^{i+1}}(\x^{i+1,*})] \leq \mathbb{E}[\langle \nabla f^{\beta^{i+1}}(\x^{i+1,0}), \x^{i+1, 0}  -\x^{i+1, *} \rangle]\\
        \leq \mathbb{E}[||\nabla f^{\beta^{i+1}}(\x^{i+1,0})|| \, ||\x^{i+1, 0}  -\x^{i+1, *}|| ].
    \end{align*}
    The result follows using the Lemma \ref{bound_grad} and since $\x^{i+1,0}$ belongs to the ball  $\mathcal{B}^{\epsilon}(\x^{i,*})$.
\end{proof}
Moreover, an estimate on the number of iterations required to reduce the norm of the gradient below some threshold may be given.
\begin{Lem}
    Under Assumptions \ref{assum1}, \ref{assum2}  and \ref{assum3}, for a subproblem $i > I$ and in the setting of Algorithm \ref{algo2}, let $s_2^{i,0} \in \R^+$ be such that $k = 1$ in Equations (\ref{cond1}) and (\ref{cond2}), assume that $L = \max( L_0(F), L_1(F))$, $\alpha_1 = \frac{3}{4}$ and $\alpha_2 = \frac{1}{2}$. Then, for a uniformly randomly chosen $R \in [0, K_i]$, it follows that
    \begin{equation*}
        \mathbb{P} \left( ||\nabla f^{\beta^i}(\x^{i,R})|| \geq \frac{L \beta^i}{4 \beta^0}  \right) \leq \frac{4 \beta^0}{\beta^i K_i^{\frac{1}{4}}} (A^i + B^i),
    \end{equation*}
    where $A^i$ and $B^i$ are defined in Equation (\ref{const_i}).
    \label{bound_prob1}
\end{Lem}
\begin{proof}
    Markov's inequality implies that
    \begin{align*}
        \mathbb{P} \left( ||\nabla f^{\beta^i}(\x^{i,R})|| \geq \frac{L \beta^i}{4 \beta^0}  \right) \leq \frac{4 \beta^0 \mathbb{E}[||\nabla f^{\beta^i}(\x^{i,R})||] }{L \beta^i}.
    \end{align*}
    Now, given the result of Theorem \ref{conv_expect} with the specified value of $\alpha_1$ and $\alpha_2$ and the fact that $L_1(f^{\beta^i})~\leq~L_1(F)$ together with Lemma \ref{lem_convex}, it follows that 
    \begin{align*}
        \frac{4 \beta^0 \mathbb{E}[||\nabla f^{\beta^i}(\x^{i,R})||] }{L \beta^i} \leq \frac{4 \beta^0}{\beta^i K_i^{\frac{1}{4}}} (A^i + B^i),
    \end{align*}
    where 
    \begin{align}
    \begin{split}
        A^i &= \frac{\rho}{s_1^{i,0}}\left( \frac{\beta^{i-1}}{4 \beta^0} + 10n \frac{s_1^{i-1,0}}{s_2^{i-1,0} K_{i-1}^{\frac{1}{4}}} + (n+3)^{\frac{3}{2}} (\beta^i - \beta^{i+1})  \right)\\
        B^i &= \frac{n s_1^{i,0}}{2} H_k^{(-\frac{3}{2})} + \ln(K_i) \left(6 \sqrt{s_2^{i,0}} (n+4) \sqrt{n} + \frac{20 n s_1^{i,0}}{s_2^{i,0}} \right),
        \end{split}
        \label{const_i}
    \end{align}
    and $K_i$ is the iteration number for subproblem $i$ and $H_k^{(-\frac{3}{2})}$ is the generalized harmonic number. 
\end{proof}
The following Lemma provides an estimate on the number of iterations required to bound the norm of the difference between $\m$ and the gradient below a certain threshold.
\begin{Lem}
\label{bound_prob2}
     For a subproblem $i \in \N$ and in the setting of Algorithm \ref{algo2}, let $s_2^{i,0} \in \R^+$ be such that $k = 1$ in Equations (\ref{cond1}) and (\ref{cond2}), assume $L = \max(L_0(F), L_1(F))$, $\alpha_1 = \frac{3}{4}$ and $\alpha_2 = \frac{1}{2}$. Then, for a uniformly randomly chosen $R \in [0, K_i]$, it follows that
    \begin{equation*}
        \mathbb{P} \left( ||\m^{i,R} - \nabla f^{\beta^i}(\x^{i,R})|| \geq \frac{L \beta^i}{4 \beta^0}  \right) \leq \frac{4 \beta^0}{\beta^i K_i^{\frac{1}{4}}} \left(3 \sqrt{s_2^{i,0}} (n +4)\sqrt{n} + \frac{10 n s_1^{i,0}}{s_2^{i,0}} \right).
    \end{equation*}
\end{Lem}
\begin{proof}
    By Markov's inequality, it follows that
    \begin{align*}
        \mathbb{P} \left( ||\m^{i,R} - \nabla f^{\beta^i}(\x^{i,R})|| \geq \frac{L \beta^i}{4 \beta^0}  \right) &\leq \frac{4 \beta^0 \mathbb{E}[||\m^{i,R} - \nabla f^{\beta^i}(\x^{i,R})||] }{L \beta^i}\\
        &= \frac{4 \beta^0}{L \beta^i K_i} \sum_{k = 0}^{K_i} \mathbb{E}[||\m^{i,k} - \nabla f^{\beta^i}(\x^{i,k})||]\\
        &\leq \frac{4 \beta^0}{ \beta^i K_i^{\frac{1}{4}}}\left(3 \sqrt{s_2^{i,0}} (n +4)\sqrt{n} + \frac{10n s_1^{i,0}}{s_2^{i,0}} \right),
    \end{align*}
    where the last inequality holds by Proposition \ref{prop2} and \ref{prop3} with $\alpha_1 = \frac{3}{4}$ and $\alpha_2 = \frac{1}{2}$.
\end{proof}

Finally, the main theorem of this section may be stated. 
\begin{Th}
\label{main_th}
    Let Assumptions \ref{assum1}, \ref{assum2} and \ref{assum3} hold and let $I$  be the threshold from Assumption \ref{assum3}.
    
    (i) For $i \in \N$,  set 
    \begin{equation*}
       \beta^i = \frac{1}{\sqrt{n}(i+1)^2}, s_1^{i,0} = \frac{1}{6 n (i+1)^{3/2}} \text{ and }   s_2^{i,0} = \frac{s_2}{(i+1)}
    \end{equation*}
    with $s_2$ so that Equations (\ref{cond1}) and (\ref{cond2}) are satisfied for $k = 1$. Moreover, let denote $K_i$ the first iteration for which $||\m^{i,K_i}|| \leq \frac{L \beta^i}{4 \beta^0}$ and that without loss of generality $L = \max(L_0(F), L_1(F))$. Let $\epsilon >0$ be a desired accuracy and let $i^* \geq \sqrt{\frac{L}{\epsilon}} \geq I$. If for any $i \geq I, K_i \geq (i+1)^6$, then after at most 
    \begin{equation*}
        O \left( \frac{n^6 L^{7/2}}{\epsilon^{7/2}} \right)
    \end{equation*}
    function evaluations, the following inequality holds
    \begin{equation}
    \label{main_equation}
        ||\nabla f^{\beta^{i^*}}(x^{i^*, 0})|| \leq \epsilon.
    \end{equation}
    
    (ii) Furthermore, when for every $i \in \N, f^{\beta^i}$ is convex then under the same setting that Theorem \ref{th_convex} given in Equation (\ref{setting}), it follows that after at most 
    \begin{equation*}
        O\left( \frac{n^\frac{9}{2} L^{7/2}}{\epsilon^{7/2}} \right)
    \end{equation*}
    function evaluations, the inequality in Equation (\ref{main_equation})  holds.
\end{Th}
\begin{proof}
    For a subproblem $i \in \N$, a probabilistic upper bound on the iteration  $ K_i \in \N$ such that $||\m^{i,K_i}|| \leq \frac{L \beta^i}{4 \beta^0}$ may be provided. We have 
    \begin{align}
    \begin{split}
        || \m^{i,K_i} ||&= \min_{k \in [0, K_i]} ||\m^{i,k}|| \\
        & \leq ||\m^{i,R}|| \\
        &\leq ||\m^{i,R} - \nabla f^{\beta^i}(\x^{i,R})||  + ||\nabla f^{\beta^i}(\x^{i,R})||,
    \end{split}
    \label{equa_min}
    \end{align}
where $R \sim \mathcal{U}[0, K_i]$. Now, probabilistic upper bounds on the number $K_i$ required to obtain that both terms in the right-hand side of the previous inequality are below $\frac{L \beta^i}{4 \beta^0}$. For the first term of the right-hand side in Equation (\ref{equa_min}), using the specified value of $s_1^{i,0}$, $s_2^{i,0}$ and $\beta^i$, Lemma \ref{bound_prob2} ensures that 
    \begin{align*}
         \mathbb{P} \left( ||\m^{i,R} - \nabla f^{\beta^i}(\x^{i,R})|| \geq \frac{L \beta^i}{4 \beta^0}  \right) &\leq\frac{4 \beta^0}{\beta^i K_i^{\frac{1}{4}}} \left(3 \sqrt{s_2^{i,0}} (n+4) \sqrt{n} + \frac{10n s_1^{i,0}}{s_2^{i,0}} \right)\\
         &\leq O \left(\frac{n \sqrt{n} (i+1)^{\frac{3}{2}}}{K_i^{\frac{1}{4}}} \right).
    \end{align*}
The second term of the right-hand side in Equation (\ref{equa_min}) depends on the value of $I$. For subproblems $i \leq I$, it follows by Markov's inequality and Theorem \ref{conv_expect} that 
\begin{align*}
    \mathbb{P} \left( ||\nabla f^{\beta^i}(\x^{i,R})|| \geq \frac{L \beta^i}{4 \beta^0}  \right) &\leq  \frac{4 \beta^0}{L \beta^i} \mathbb{E}[||\nabla f^{\beta^i}(\x^{i,R})||]\\
    &\leq \frac{4 \beta^0}{ \beta^i} \left( \frac{D_f^i}{s_1^{i,0}} + \frac{n s_1^{i,0}}{2}H_k^{{(-\frac{3}{2})}} + \ln(K_i) \left( 6  \sqrt{s_2^{i,0}} (n +4) \sqrt{n} + \frac{40  s_1^{i,0} n}{s_2^{i,0}}\right) \right)\\
    &\leq O\left(  \frac{\max \left( \frac{n(i+1)^{\frac{7}{2}}}{L}, n \sqrt{n} \ln(K_i)(i+1)^{\frac{3}{2}}\right)}{K_i^{\frac{1}{4}}} \right).
\end{align*}
For subproblems $i > I$, Lemma \ref{bound_prob1}  ensures that
\begin{align*}
    \mathbb{P} \left( ||\nabla f^{\beta^i}(\x^{i,R})|| \geq \frac{L \beta^i}{4 \beta^0}  \right) \leq   \frac{4 \beta^0}{\beta^i K^{\frac{1}{4}}} (A^i + B^i),
\end{align*}
where $A^i$ and $B^i$ are given in Equation (\ref{const_i}). Now, given condition on $K_i$, it follows that
    \begin{align*}
        A^i &= \rho n (i+1)^{3/2} \left( \frac{1}{i^2} + \frac{10}{s_2 i^2} + \frac{2 (n+3)}{i^2(i+1)} \right) \text{ and }\\
        B^i &= \frac{H_k^{(-\frac{3}{2})}}{2(i+1)^{3/2}}    + \ln(K_i) \left( \frac{6n \sqrt{n+3} \sqrt{s_2}}{\sqrt{i+1}} + \frac{12}{s_2 \sqrt{i+1} }\right).
    \end{align*}
    Thus, we obtain
    \begin{equation}
        \mathbb{P} \left( ||\nabla f^{\beta^i}(\x^{i,R})|| \geq \frac{L \beta^i}{4 \beta^0}  \right) \leq O \left( \frac{ n \sqrt{n} (i+1)^{\frac{3}{2}} \ln(K_i) }{K_i^{\frac{1}{4}}} \right).
    \end{equation}
    Therefore,  to obtain $||\m^{i,K_i}|| \leq \frac{L \beta^i}{4 \beta^0}$, it takes at most 
    \begin{equation*}
        K_i = \left\{
        \begin{array}{ll}
            O \left(\max \left( n^4 (i+1)^{14}, n^6 (i+1)^6 \right) \right) & \mbox{if  } i \leq I \\
             O \left( \left(n^6 (i+1)^6 \right) \right) & \mbox{otherwise,}
        \end{array}
\right.
    \end{equation*}
    iterations. Thus, by taking $i^* \geq \sqrt{\frac{L}{\epsilon}}$, it follows that the number of iterations needed to reach the subproblem $i^*$ is
    \begin{align}
    \begin{split}
           \label{equa_O}
        \sum_{i = 1}^{i^*} K_i &= \sum_{i = 1}^{I } K_i  + \sum_{i = I+1}^{i^*} K_i \\ &=  O \left(\max \left( n^4 (I+1)^{15}, n^6 (I+1)^7 \right) \right) +  O(n^6 (i^*)^7) \\ &= O \left( \frac{n^6 L^{7/2}}{\epsilon^{7/2}}\right),
    \end{split}
    \end{align}
    where $I $ is a constant with respect to $\epsilon$. Once this number of iterations is reached, it follows that  $||\m^{i^*,0}|| \leq \frac{L}{(i^*+1)^2} \leq \epsilon$ and  by Lemma \ref{bound_grad}
    \begin{equation*}
        ||\nabla f^{\beta^{i^*}}(\x^{i^*, K_{i^*}})|| \leq \frac{L}{(i^*+1)^2} + \frac{L}{\sqrt{i^* +1} (i^*)^\frac{3}{2}} \leq 2 \epsilon.
    \end{equation*}
    For the second part of the proof, the bounds on Equation (\ref{equa_min}) does not depend on the value of $I$ since $f^{\beta^i}$ is assumed convex for every $i \in \N$. With the setting in Equation (\ref{setting}), it follows that 
    \begin{align*}
        \mathbb{P}\left(||\nabla f^{\beta^i}(\x^{i,R})|| \geq \frac{L \beta^i}{4 \beta^0} \right) &\leq \frac{4 \beta^0}{\beta^i} \mathbb{E}[||\nabla f^{\beta^i}(\x^{i,R}) ||] \leq 16 \frac{\rho n \sqrt{n} (i+1)^2}{K_i^{\frac{1}{3}}} \text{ and }\\
        \mathbb{P} \left( ||\m^{i,R} - \nabla f^{\beta^i}(\x^{i,R})|| \geq \frac{L \beta^i}{4 \beta^0}  \right) &\leq \frac{4 \beta^0}{L \beta^i} n\sqrt{n} L \frac{ \sum_{k = 0}^{K_i-1} \frac{2 \rho}{\Gamma^{k+1}(k+1)^{\frac{4}{3}}}}{\sum_{k =0}^{K_i-1} \frac{2 \rho}{\Gamma^{k+1}(k+1)}} \leq 8 \frac{n \sqrt{n} (i+1)^2}{K_i^{\frac{1}{3}}},
    \end{align*}
    where the first inequality follows by Theorem \ref{th_convex} and the second one by the definition of the probability density of $R$ together with Propositions \ref{prop2} and \ref{prop3}. Therefore, it takes at most $K_i = O(n^{\frac{9}{2}}(i+1)^6 )$ to obtain $||\m^{i,K_i}|| \leq \frac{L \beta^i}{4 \beta^0}$. Thus, by taking $i^* \geq \sqrt{\frac{L}{\epsilon}}$, it follows that the number of iterations needed to reach the subproblem $i^*$ is 
    \begin{equation*}
        \sum_{i = 1}^{i^*} K_i = O(n^{\frac{9}{2}}(i^*)^7) = O \left(\frac{n^{\frac{9}{2}} L^{\frac{7}{2}}}{\epsilon^{\frac{7}{2}}} \right).
    \end{equation*}
  It remains to apply the Lemma \ref{bound_grad} as previously to complete the proof.
\end{proof}

We would like to make a few remarks about this theorem. 
First, one approach to satisfy the condition $K_i \geq (i+1)^6$ for any $i \in \N$ is to incorporate it into the stopping criterion of Algorithm \ref{algo1}. 
However, due to the limited number of iterations in practice, this condition is typically replaced by a weaker one, $ K_i \geq M$, where $M > 0$ is a constant. 
Second, the main result of Theorem \ref{bound_grad} establishes an $\epsilon$ convergence rate for a single run of the SSO algorithm, which is the first of its kind to the best of our knowledge. 
This was made possible by decomposing the problem given in Equation~(\ref{prob1}) into a sequence of subproblems, each of which is solved using carefully chosen stopping criteria and step sizes. It is worth noting that, in \cite{ghadimi2013stochastic}, the $(\epsilon, \Lambda)$-solution of the norm of the gradient is obtained after at most $O\left(\frac{n L^2 \sigma^2}{\epsilon^4}\right)$. 
Although this bound has a weaker dependence on $n$ and $L$, it is worse in terms of $\epsilon$. 
Third, the first term in Equation (\ref{equa_O}) may be significant even if it is fixed, particularly if the region where the function is convex is difficult to reach, indeed this constant disappears when $f^{\beta^i}$ is convex for every index $i$. Nevertheless, the bounds given are the worst one and may be considerably smaller in practice. 
A way to decrease this term is to decrease the power on $i$ in the denominator of $\beta^i$, $s_1^{i,0}$ and $s_2^{i,0}$ but it also decreases the asymptotic rate of convergence. 
Finally, the process used in the SSO algorithm may be extended to other momentum-based methods and give an appealing property for these methods compared to the classical SGD.  

\section{Numerical experiments}
\label{sec-num_res}

The numerical experiments are conducted for two bounded constrained blackbox optimization problems. In order to handle the bound constraints $\x \in [\ell, \u] \subset \R^n$, the update in Equation (\ref{update_x}) is simply projected such that 
$\x \leftarrow \max(\ell, \min(\x, \u))$.

\subsection{Application to a solar thermal power plant}
The first stochastic test problem is SOLAR \footnote{https://github.com/bbopt/solar} \cite{MScMLG}, which simulates a thermal solar power plant and contains several instances allowing to choose the number of variables, the types of constraints and the objective function to optimize. 
All the instances of SOLAR are stochastic, have nonconvex constraints and integer variables. 
In this work, the algorithms developed does not deal with integer variables. 
Therefore, the problem is altered: all integer variables are fixed to their initial value and the problem is to obtain a feasible solution by optimizing the expectation of constraint violations over the remaining variables. 
Numerical experiments are conducted for the second instance of the SOLAR framework, which considers 12 variables (2 integers) and 12 constraints:
\begin{equation*}
    \min_{\x \in [0,1]^{12}} \mathbb{E} \left[ \sum_{j = 1}^{m} \max(0, c_j(\x, \boldsymbol{\xi}))^2 \right]
\end{equation*}
where the $c_j$ are the original stochastic constraints and the bound constraints have been normalized. The second instance of SOLAR is computationally expensive; a run  may take between several seconds and several minutes. Therefore, the maximum number  of function evaluations is set to $1000$. Four algorithms are used:
\begin{itemize}
    \item SSO, whose the hyperparameters values are given in Table \ref{tab4}. The search step given in Algorithm \ref{algo2} is used for this experiment. A truncated version of the Gaussian gradient based estimate is used for this experiment.
    \item ZO-adaMM \cite{chen2019zo} which is a zeroth-order version of the original Adam algorithm. This algorithm appears as one of the most effective according to \cite{liu2020primer, chen2019zo} in terms of distortion value, number of function evaluations and success rate. The default parameters defined experimentally in \cite{chen2019zo} are used on this problem, except that $\beta = 0.05$ and the learning rate is equal to $0.3$. Moreover, the same gradient estimator that ZO-Signum is used to eliminate its impact on the performance. 
    \item CMA-ES \cite{Hansen2006} an algorithm based on biological inspired operators. Its name comes from the adaptation of the covariance matrix of the multivariate normal distribution used during the mutation. The version of CMA-ES used is the one of the pymoo \cite{pymoo} library with the default setting.
    \item The NOMAD 3.9.1 software \cite{Le09b}, based on the Mesh Adaptive Direct Search (MADS) \cite{AuDe2006} algorithm, a popular blackbox optimization solver. 
\end{itemize}

\begin{figure}[ht!]
    \centering
    \includegraphics[width = 0.6 \linewidth]{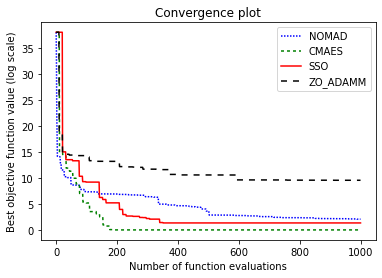}
    \caption{Average of 5 different seed runs for the NOMAD, CMAES, SSO and ZO-adaMM algorithms.}
    \label{convergence_plot}
\end{figure}
The results are presented in Figure \ref{convergence_plot}, which plots the average best result obtained by each algorithm with five different seeds. In this experiment, SSO obtains similar performance to NOMAD and CMAES which are state-of-the-art algorithms for this type of problem.  ZO-adaMM has difficulty to converge even though it is  a ZO algorithm. 

\subsection{Application to blackbox adversarial attack}

This section demonstrates the competitiveness of the SSO algorithm through experiments involving the generation of blackbox adversarial examples for Deep Neural Networks (DNNs) \cite{xu2018structured}. Generating an adversarial example for a DNN involves adding a well-designed perturbation to the original legal input to cause the DNN to misclassify it. In this work, the attacker considers the DNN model to be unknown, hence the term blackbox. Adversarial attacks against DNNs are not just theoretical, they pose a real safety issue \cite{papernot2017practical}. Having an algorithm that generates effective adversarial examples enables modification of DNN architecture to enhance its robustness against such attacks. An ideal adversarial example is one that can mislead a DNN to recognize it as any target image label, while appearing visually similar to the original input, making the perturbations indiscernible to human eyes. The similarity between the two inputs is typically measured by a $\ell_p$ norm. Mathematically, a blackbox adversarial attack can be formalized as follows. Let $(\y, \ell)$ denote a legitimate image $\y$ with the true label $\ell \in [1, M]$, where $M$ is the total number of image classes. Let $\x$ denote the adversarial perturbation; the adversarial example is then given by $\y' = \y + \x$, and the goal is to solve the  problem \cite{chen2019zo}
\begin{align*}
    &\min_{\x} \lambda f(\y + \x) + ||\x||_2 \\
    &\text{ subject to } (\y + \x) \in [-0.5, 0.5]^n,
\end{align*}
where $\lambda > 0$ is a regularization parameter and $f $ is the blackbox attack loss function. In our experiments, $\lambda = 10$ and the loss function is defined for untargeted attack \cite{carlini2017towards}, i.e,
\begin{equation*}
    f(\y') = \max \{ Z(\y')_\ell - \max_{j \neq \ell} Z(\y)_j, 0\},
\end{equation*}
where $Z(\y')_k$ denotes the prediction score of class $k$ given the input $\y'$. Thus, the minimum value of $0$ is reached as the perturbation succeeds to fool the neural network. 

The experiments of generating blackbox adversarial examples are first performed on an adapted AlexNet \cite{krizhevsky2017imagenet} under the dataset Cifar10 and then on InceptionV3 \cite{szegedy2016rethinking} under the dataset ImageNet \cite{deng2009imagenet}. 
Since the NOMAD algorithm is not recommended for large problems, three algorithms are compared  : SSO (without search),  ZO-adaMM and CMAES. 
In the experiments, the hyperparameters of the algorithm ZO-adaMM are taken as in \cite{chen2019zo}, those of SSO are given in Table \ref{tab4} and the uniform gradient based estimate is used for both algorithms. 
Moreover, for the Cifar10 dataset, different initial learning rates for ZO-adaMM are used to observe its influence on the success rate. 
Experiments are conducted for 100 randomly selected images with a starting point corresponding to a null distortion, the maximum number of function queries is set to $5000$. Thus, as the iteration increases, the attack loss decreases until it converges to $0$ (indicating a successful attack) while the norm of the distortion could increase. 

The best attack performance involves a trade-off between a fast convergence to a $0$ attack loss in terms of function evaluations, a high rate of success, and a low distortion (evaluated by the $\ell_2$-norm). 
The results for the Cifar10 dataset are given in Table \ref{tabcifar}. 
\begin{center}
\small
    \captionof{table}{Results of blackbox adversarial attack for the Cifar10 dataset ($n= 3\times32\times32$) }
    \begin{tabular}{|C{3cm}|C{3cm}|C{3cm}|C{3cm}|C{3cm}|}
    \hline  Method & Attack success rate & $||\ell_2||$ first success&   Average \# of function evaluations  \\
    \hline
    ZO-adaMM $lr = 0.01$& 79 \% & $0.14$ &$582$\\
    \hline
    ZO-adaMM $lr = 0.03$& 96\% & $0.97$ & $310$\\
    \hline
    ZO-adaMM $lr = 0.05$& 98\% & $2.10$ & $215$\\
    \hline
    CMAES $\sigma = 0.005$& 99\% & $0.33$ &  $862$\\
    \hline
    SSO & 100\% & $0.55$ &  $442$\\
    \hline
    \end{tabular} 
    \label{tabcifar}
\end{center}

Except for ZO-adaMM with an initial learning rate equal to $0.01$, all algorithms achieve a success rate above $95 \%$. Among these algorithms, ZO-adaMM with a learning rate equal to $0.05$, has the best convergence rate in terms of function evaluations but has the worst value of distortion. On the contrary, CMA-ES obtains the best value of distortion but has the worst convergence rate. The SSO algorithm obtains balanced results, and is the only one that reaches full success rate. 

Table \ref{tabimagenet} displays results for the ImageNet dataset. 
Only two algorithms are compared since the dimension is too large to invert the covariance matrix in CMA-ES. For this dataset, ZO-adaMM and SSO have the same  convergence rate. However, SSO outperforms ZO-adaMM in terms of success rate while having a slightly higher level of distortion.
\begin{center}
\small
    \captionof{table}{Results of blackbox adversarial attack for the ImageNet dataset ($n = 3 \times 299 \times 299$) }
    \begin{tabular}{|C{3cm}|C{3cm}|C{3cm}|C{3cm}|C{3cm}|}
    \hline  Method & Attack success rate & $||\ell_2||$ first success&   Average \# of function evaluations  \\
    \hline
    ZO-adaMM $lr = 0.01$& $59$ \% & $19$ &$1339$\\
    \hline
    SSO & $73$ \% & $33$ &  $1335$\\
    \hline
    \end{tabular} 
    \label{tabimagenet}
\end{center}

\section{Concluding remarks}
\label{sec-conclusion}

This paper presents a method for computationally expensive stochastic blackbox optimization. The approach uses zeroth-order gradient estimates, which provides three advantages. 
First, they require few function evaluations to estimate the gradient, regardless of the problem's dimension. 
Second, under mild conditions on the noised objective function, the problem is formulated as optimizing a smoothed functional. 
Third, the smoothed functional may appear to be locally convexified near a local minima.

Based on these three features, the SSO algorithm was proposed. This algorithm is a sequential one and comprises two steps. The first is an optional search step that improves the exploration of the decision variable space and the algorithm's efficiency. The second is a local search, which ensures the convergence of the algorithm. In this step, the original problem is decomposed into subproblems solved by a ZO-version of a sign stochastic descent with momentum algorithm. 
More specifically, when the momentum's norm falls below a specified threshold that depends on the smoothing parameter, the subproblem is considered solved. The smoothing parameter's value is then decreased, and the SSO algorithm moves on to the next subproblem.

A theoretical analysis of the algorithm is conducted. 
Under Lipschitz continuity of the stochastic zeroth-order oracle, a convergence rate in expectation of the ZO-Signum algorithm is derived. 
Under additional assumptions of smoothness and  convexity or local convexity of the objective function near its minima, a convergence rate of the SSO algorithm to an $\epsilon$-optimal point of the problem is derived, which is, to the best of our knowledge, the first of its kind.

Finally, numerical experiments are conducted on a solar power plant simulation and on adversarial blackbox attacks. Both examples are computationally expensive,  the former is a small size problem ($n \approx 10$) and the latter is a large size problem (up to $n \approx 10^5$). 
The results demonstrate the SSO algorithm's competitiveness in both performance and convergence rate compared to state-of-the-art algorithms. Further work will extend this approach to constrained stochastic optimization.






\newpage
\appendix

\section{Notations}
The following list describes symbols used within the body of the document. Throughout the paper, when a symbol is shown in bold then it is a vector, otherwise it is a scalar.
\footnotesize
\makenomenclature
\nomenclature[01]{$n$}{The dimension of the space of the design variables}
\nomenclature[01]{$\Omega$}{The sample space of $\boldsymbol{\xi}$, i.e, the set of all possible outcomes of $\boldsymbol{\xi}$}.
\nomenclature[02]{$\boldsymbol{\xi}: \Omega \to  \R^m $}{The vector of uncertainties}
\nomenclature[03]{$\mathbb{E}_{\boldsymbol{\xi}}[\cdot]$}{The expectation with respect to the random vector $\boldsymbol{\xi}$}
\nomenclature[03]{$F : \R^n\times \R^m \to \R$}{The stochastic zeroth-order oracle that takes into account the uncertainty $\boldsymbol{\xi}$}
\nomenclature[04]{$f : \R^n \to \R$}{The expectation of $F$ with respect to $\boldsymbol{\xi}$}
\nomenclature[05]{$\beta \in \R^{+*}$}{A strictly positive scalar using as smoothing parameter}
\nomenclature[06]{$\u \in \R^n$}{A Gaussian random vector}
\nomenclature[07]{$f^\beta = \mathbb{E}[f(\x + \beta \u)]$}{A smooth approximation of a function $f$}
\nomenclature[08]{$L_0(f)$}{The Lipschitz constant associated to a function $f$}
\nomenclature[09]{$L_1(f)$}{The Lipschitz constant associated to the gradient of a function $f$}
\nomenclature[10]{$\nabla f $}{The gradient of a function $f$}
\nomenclature[11]{$\Tilde{\nabla} f $}{An estimator of the gradient of a function $f$}
\nomenclature[11]{$\Tilde{\g} $}{An estimator of the gradient of a function $f$ based on outputs of the stochastic zeroth-order oracle $F(\x, \boldsymbol{\xi})$}
\nomenclature[12]{$j \in [1, n]$}{The counter associated with the dimension}
\nomenclature[12]{$i \in \N $}{The outer iteration counter associated with a subproblem}
\nomenclature[13]{$k \in \N $}{The inner iteration counter}
\nomenclature[14]{$\m \in \R^n $}{The momentum vector}
\nomenclature[16]{$s_2^{i,k} \in (0, 1)$}{The step size associated with the momentum }
\nomenclature[17]{$s_1^{i,k} \in (0, 1)$}{The step size associated with $\x$}
\nomenclature[18]{$L \in \R^{+*}$}{An approximation of the Lispchitz constant }
\nomenclature[19]{$q \in \N$}{The size of the mini batch used to estimate $\Tilde{\nabla}$ }
\nomenclature[20]{$M \in \N$}{The minimum number of iteration used in the ZOS algorithm}
\nomenclature[21]{$H_k^{(\alpha)} $}{The generalized harmonic number of order $\alpha$}
\nomenclature[22]{$\mathcal{C}^{0+}$}{Class of Lipschitz continuous functions}
\nomenclature[23]{$\mathcal{C}^{1+}$}{Class of differentiable functions whose the gradient is Lipschitz}
\nomenclature[23]{$\mathcal{C}^{\infty}$}{Class of infinitely differentiable functions}
\printnomenclature[3.1 cm]
\normalsize
\newpage
\section{Proof of Proposition \ref{prop1}}
\label{appendix_A}
\begin{Prop}[\cite{bernstein2018signsgd}]
For the subproblem $i \in \N$, under Assumption \ref{assum1} and in the setting of Algorithm \ref{algo1}, we have
\begin{align}
\begin{split}
        s_1^{i,k} \mathbb{E}[||\nabla f^{\beta^i} (\x^{i,k}) ||_1] \leq &\mathbb{E}[f^{\beta^i}(\x^{i,k}) - f^{\beta^i}(\x^{i,k+1})]  + \frac{n L_1(f^{\beta^i})}{2} (s_1^{i,k})^2 \\
    &+ 2 s_1^{i,k} \underbrace{\mathbb{E}[||\Bar{\m}^{i,k+1}  - \nabla f^{\beta^i}(\x^{i,k})||_1]}_{\text{bias}} + 2s_1^{i,k}  \sqrt{n} \sqrt{\underbrace{\mathbb{E}[||\m^{i,k+1} - \Bar{\m}^{i,k+1}||_2^2]}_{\text{variance}}} 
\end{split}
\end{align}
where $\Bar{m}_j^{i,k+1}$ is defined recursively as $\Bar{m}_j^{i,k+1} = s_2^{i,k} \nabla f^{\beta^i}(\x^{i,k}) + (1-s_2^{i,k})\Bar{m}_j^{i,k} $.
\end{Prop}
\begin{proof}
    By $L_1(f^{\beta^i})$-Lipschitz smoothness of $f^{\beta^i}$ (see Lemma \ref{lem_trunc_gaus}.3), it follows that 
    \begin{align*}
        f^{\beta^i}(\x^{i, k+1} ) &\leq f^{\beta^i}(\x^{i,k}) + \langle \nabla f^{\beta^i}(\x^{i,k}), \x^{i, k+1} - \x^{i,k} \rangle + \frac{L_1(f^{\beta^i})}{2} ||\x^{i, k+1}  - \x^{i,k} ||_2^2 \\
        &= f^{\beta^i}(\x^{i,k}) - s_1^{i,k} \langle \nabla f^{\beta^i}(\x^{i,k}), \text{sign}(\m^{i,k+1})\rangle + \frac{L_1(f^{\beta^i}) (s_1^{i,k})^2}{2} ||\text{sign}(\m^{i,k+1})||_2^2\\
        &= f^{\beta^i}(\x^{i,k}) - s_1^{i,k} ||\nabla f^{\beta^i}(\x^{i,k}) ||_1 + \frac{n L_1(f^{\beta^i}) }{2}(s_1^{i,k})^2 \\
        &+ 2s_1^{i,k} \sum_{j = 1}^n |\nabla_j f^{\beta^i}(\x^{i,k})| \mathbf{1} \{\text{sign}(m_j^{i,k+1}) \neq \text{sign}(\nabla_j f^{\beta^i}(\x^{i,k}) )\},
    \end{align*}
    where $\mathbf{1}\{\cdot\}$ is the indicator function. Now, as in \cite{bernstein2018signsgd, liu2018signsgd}, the expected improvement conditioned on $\x^{i,k}$ is given by
    \begin{align}
    \mathbb{E}[f^{\beta^i}(\x^{i,k+1}) - f^{\beta^i}(\x^{i,k})| \x^{i,k}] \leq &-s_1^{i,k} ||\nabla f^{\beta^i} (\x^{i,k}) ||_1  + \frac{n L_1(f^{\beta^i})}{2} (s_1^{i,k})^2 \\
    &+ 2s_1^{i,k} \sum_{j = 1}^n |\nabla_j f^{\beta^i}(\x^{i,k})| \mathbb{E}[\mathbf{1}\{\text{sign}(m_j^{i,k+1}) \neq \text{sign}(\nabla_j f^{\beta^i}(\x^{i,k}) ) \} | \x^{i,k}].
    \label{first_step}
    \end{align}
    Again, as in \cite{bernstein2018signsgd, liu2018signsgd}, the expectation that the sign of $m_j^{i,k+1}$ be different of the sign of $\nabla_j f^{\beta^i}(\x^{i,k})$ is relaxed by considering that the set
    \begin{equation*}
        \{ m_j^{i,k+1} : \text{sign}(m_j^{i,k+1}) \neq \text{sign}(\nabla_j f^{\beta^i}(\x^{i,k}) \} \subset \{m_j^{i,k+1} : \; |m_j^{i,k+1} - \nabla_j f^{\beta^i}(\x^{i,k})| \geq |\nabla_j f^{\beta^i}(\x^{i,k})| \}.
    \end{equation*}
    Therefore, it follows that
    \begin{align}
        \mathbb{E}[\mathbf{1}\{\text{sign}(m_j^{i,k+1}) \neq \text{sign}(\nabla_j f^{\beta^i}(\x^{i,k}) )\} |\x^{i,k}] &\leq \mathbb{E}[\mathbf{1} \{ |m_j^{i,k+1} - \nabla_j f^{\beta^i}(\x^{i,k})| \geq |\nabla_j f^{\beta^i}(\x^{i,k})|\} |\x^{i,k}] \\
        &\leq \frac{\mathbb{E}[|m_j^{i,k+1} - \nabla_j f^{\beta^i}(\x^{i,k})|\; |\x^{i,k}]}{ |\nabla_j f^{\beta^i}(\x^{i,k})|},
        \label{markov}
    \end{align}
    where the second inequality comes from conditional Markov's inequality. Substituting Equation (\ref{markov}) into Equation (\ref{first_step}) and taking expectation over all the randomness we obtain
    \begin{align}
    \begin{split}
         \mathbb{E}[f^{\beta^i}(\x^{i,k+1}) - f^{\beta^i}(\x^{i,k})] \leq &-s_1^{i,k} \mathbb{E}[||\nabla f^{\beta^i} (\x^{i,k}) ||_1 ]  + \frac{nL}{2} (s_1^{i,k})^2 \\ &+ 2s_1^{i,k} \sum_{j = 1}^n  \mathbb{E}[|m_j^{i,k+1} - \nabla_j f^{\beta^i}(\x^{i,k})| ].
    \end{split}
    \label{step_2}
    \end{align}
    Moreover, by adding and subtracting $\Bar{\m}^{i,k+1}$ in the terms of the sum of Equation (\ref{step_2}),
    \begin{align*}
         \sum_{j = 1}^n  \mathbb{E}[|m_j^{i,k+1} - \nabla_j f^{\beta^i}(\x^{i,k})|] &=    \mathbb{E}[||\m^{i,k+1} - \Bar{\m}^{i,k+1} + \Bar{\m}^{i,k+1} -  \nabla f^{\beta^i}(\x^{i,k})||_1]\\
         &\leq  \sqrt{n} \mathbb{E}[||\m^{i,k+1} - \Bar{\m}^{i,k+1} ||_2] +  \mathbb{E}[||\Bar{\m}^{i,k+1} -  \nabla f^{\beta^i}(\x^{i,k})||_1]\\
         &\leq \sqrt{n}\sqrt{\mathbb{E}[||\m^{i,k+1} - \Bar{\m}^{i,k+1}||_2^2]} + \mathbb{E}[||\Bar{\m}^{i,k+1} -  \nabla f^{\beta^i}(\x^{i,k})||_1],
    \end{align*}
    where the first inequality comes from $||\cdot||_1 \leq \sqrt{n}||\cdot||_2$ and the second one from Jensen's inequality. Finally, incorporating  the last inequality in Equation (\ref{step_2}) completes the proof.
\end{proof}

\section{List of hyperparameters for the SSO algorithm}
\begin{center}
\small
    \captionof{table}{List of hyperparameters for the SSO algorithm}
    \begin{tabular}{|C{2.5cm}|C{2.5cm}|C{2.5cm}|C{2.5cm}|C{2.5cm}|C{2.5cm}|}
    \hline  Problem & $\beta^i$ & $s_1^{i,k}$ & $s_2^{i,k}$ &  $M$ & $q$ \\
    \hline
     Cifar10 & $\frac{0.005}{(i+1)^2}$  & $\frac{0.005}{(i+1)^{\frac{3}{2}}\sqrt{k+1}}$ &$\frac{0.9}{(i+1) (k+1)^{\frac{1}{4}}}$ & $60$ & $10$\\
    \hline
    ImageNet &  $\frac{0.001}{(i+1)^{2}}$  & $\frac{0.003}{(i+1)^{\frac{3}{2}}\sqrt{k+1}}$ &$\frac{0.7}{(i+1) (k+1)^{\frac{1}{4}}}$ & $100$ &$10$\\
    \hline
    Solar  &  $\frac{0.3}{(i+1)^2}$  & $\frac{0.1}{(i+1)^{\frac{3}{2}}\sqrt{k+1}}$ & $\frac{0.5}{(i+1) (k+1)^{\frac{1}{4}}}$ & $5$ & $10$\\
    \hline
    \end{tabular} 
    \label{tab4}
\end{center}

\section{Original signSGD and Signum algorithms}
\label{appendix_C}

Below are the original versions of the signSGD and Signum algorithms.
\begin{algorithm}[htb!]
	\caption{signSGD algorithm} 
	\begin{algorithmic}[1]
        \State{\textbf{Input:} $\x^{0}$, $s_1 \in (0, 1)$}
        \For{ $k = 0, 1,  \dots$}
        \State{Calculate an estimate of the stochastic gradient $\Tilde{\nabla} f(\x^k)$ and update:}
        \begin{align*}
            &\x^{k+1} = \x^{k} - s_1 \text{sign}(\Tilde{\nabla} f(\x^k))  
        \end{align*}
        \EndFor
    \State{Return $\x^k$}
	\end{algorithmic} 
\end{algorithm}

\begin{algorithm}[htb!]
	\caption{Signum algorithm} 
	\begin{algorithmic}[1]
        \State{\textbf{Input:} $\x^{0}, \m^{0}, s_1 \in (0, 1), s_2 \in (0, 1)$}
        \For{ $k = 0, 1,  \dots$}
         \State{Calculate an estimate of the stochastic gradient $\Tilde{\nabla} f(\x^k)$ and update:}
        \begin{align*}
            &\m^{k+1} = s_2 \m^{k} + (1-s_2) \Tilde{\nabla} f(\x^k)  \\
            &\x^{k+1} = \x^k - s_1 \text{sign}(\m^{k+1})
        \end{align*}
        \EndFor
    \State{Return $\x^k$}
	\end{algorithmic} 
\end{algorithm}

\clearpage
\bibliographystyle{spmpsci} 
\bibliography{bibliography}
\end{document}